\documentclass[11pt,reqno]{amsart}
\usepackage{fullpage}
\usepackage[mathscr]{eucal}
\usepackage{mathrsfs}
\usepackage{amsfonts}
\usepackage{amsmath}
\usepackage{amsthm}
\usepackage{amssymb}
\usepackage{latexsym}
\usepackage[all]{xy}
\usepackage{pstricks,pst-node}
\usepackage[mathscr]{eucal}

\theoremstyle{plain}
\newtheorem{defn}[equation]{Definition}
\newtheorem{cor}[equation]{Corollary}
\newtheorem{lem}[equation]{Lemma}
\newtheorem{prop}[equation]{Proposition}
\newtheorem{thm}[equation]{Theorem}
\newtheorem{conj}[equation]{Conjecture}

\theoremstyle{remark}

\newtheorem{rem}[equation]{Remark}
\newtheorem{rems}[equation]{Remarks}
\newtheorem{examp}[equation]{Example}

\newtheorem{sublem}[equation]{Sublemma}
\numberwithin{equation}{subsection}

\makeatletter
\renewcommand{\subsection}{\@startsection{subsection}{2}{0pt}{-3ex
plus -1ex minus -0.2ex}{-2mm plus -0pt minus
-2pt}{\normalfont\bfseries}} \makeatother

\newcommand{\scr}[1]{\mathscr{#1}}

\newcommand{\dis}{\displaystyle}

\newcommand{\beq}{\begin{equation}\label}
\newcommand{\eeq}{\end{equation}}
\newcommand{\fra}[2]{\mbox{$\frac{{#1}}{{#2}}$}}

\newcommand{\hd}{{\:\raisebox{2pt}{\text{\circle*{1.5}}}}}
\newcommand{\hdot}{{\:\raisebox{2pt}{\text{\circle*{1.5}}}}}
\newcommand{\idot}{{\:\raisebox{2pt}{\text{\circle*{1.5}}}}}

\newcommand{\bimod}[1]{#1\text{-}{\sf{bimod}}}
\newcommand{\erem}{\hfill$\lozenge$\end{rem}\vskip 3pt }

\newcommand{\dok}{\vskip 2pt $\dis\begin{array}{r} \hskip }
\newcommand{\edok}{\end{array}$\qed}
\DeclareMathOperator{\Spec}{\mathrm{Spec}}

\DeclareMathOperator{\Sym}{\mathrm{Sym}}

\DeclareMathOperator{\Aut}{\mathrm{Aut}}
\DeclareMathOperator{\rk}{\mathrm{rk}}
\DeclareMathOperator{\Image}{\mathrm{Image}}
\DeclareMathOperator{\Ker}{\mathrm{Ker}}
\DeclareMathOperator{\Hom}{\mathrm{Hom}}

\DeclareMathOperator{\Perf}{\mathsf{Perf}}
\DeclareMathOperator{\Tails}{\mathsf{Tails}} \DeclareMathOperator{\rees}{{\mathsf R}}
\newcommand{\Ext}{\text{Ext}}

\newcommand{\oper}{\operatorname}

\newcommand{\vol}{{\mathsf{vol}}}
\newcommand{\ups}{\Upsilon }

\newcommand{\HH}{{H\!H}}
\def\BB{{\mathsf{B}}}
\newcommand{\PH}{{P\!H}}
\newcommand{\harr}{{\textsl{Harr}}}

\newcommand{\Id}{\text{Id}}

\newcommand{\kh}{\k[[\hbar]]}

\newcommand{\llb}{{\boldsymbol\langle}\!^{\!}{\boldsymbol\langle}}
\newcommand{\rrb}{{\boldsymbol\rangle}\!^{\!}{\boldsymbol\rangle}}
\newcommand{\xxx}{x_1,\ldots,x_n }
\def\pb{$\bullet\quad$\parbox[t]{140mm}}
\newcommand{\nnn}{1,\ldots,n }
\newcommand{\hilb}{{\mathsf P}}

\newcommand{\lan}{\langle }
\newcommand{\ran}{\rangle }
\newcommand{\pp}{\varpi}
\newcommand{\eps}{\varepsilon}
\newcommand{\mR}{{\mathcal R}}

\def\ct{\C^\times}

\newcommand{\Om}{\Omega}
\newcommand{\De}{{\Delta}}
\renewcommand{\d}{ d }
\def\ccirc{{{}_{\,{}^{^\circ}}}}
\newcommand{\Ga}{\Gamma }
\renewcommand{\div}{{\mathrm{div}}}

\newcommand{\J}{\mathbb{J}}

\renewcommand{\r}{{\mathsf R}}
\newcommand{\raa}{{\mathcal {RA}}}

\newcommand{\rkk}{{R(\!(\h)\!)}}
\newcommand{\h}{\hbar }
\newcommand{\hh}{[[\hbar]]}
\newcommand{\hb}{(\!(\hbar)\!)}

\newcommand{\wh}{\widehat }

\renewcommand{\o}{\otimes }
\newcommand{\nat}{_{\small\operatorname{cyc}}}

\newcommand{\kk}{\Bbbk}
\newcommand{\qq}{\kk\bigotimes_{\k[[\h]]}}

\newcommand{\La}{{\mathsf{\Lambda}}}
\newcommand{\iso}{{\,\stackrel{_\sim}{\to}\,}}
\newcommand{\cd}{\!\cdot\!}

\newcommand{\mm}{{\mathcal{M}}}
\newcommand{\wt}{\widetilde }

\newcommand{\ben}{\begin{enumerate}}
\newcommand{\een}{\end{enumerate}}
\newcommand{\X}{{\mathfrak{X}}}

\newcommand{\KS}{\oper{KS}}

\newcommand{\pqr}{{^{\tiny{\!P,Q,R}}}}

\newcommand{\pa}{\partial }
\newcommand{\Q}{{\mathbb{Q}}}
\newcommand{\zz}{{\mathcal Z}}

\newcommand{\E}{{\mathbb{E}}}

\newcommand{\dr}{{\mathsf{d}}}

\newcommand{\inv}{^{-1}}
\newcommand{\vi}{${\en\sf {(i)}}\;$}
\newcommand{\vii}{${\;\sf {(ii)}}\;$}
\newcommand{\viii}{${\sf {(iii)}}\;$}

\newcommand{\sset}{\subset}

\newcommand{\sminus}{\smallsetminus}
\newcommand{\mto}{\mapsto}
\newcommand{\into}{{}^{\,}\hookrightarrow^{\,}}
\newcommand{\too}{\,\longrightarrow\,}
\newcommand{\onto}{\twoheadrightarrow}

\DeclareMathOperator{\gr}{\mathrm{gr}}

\def\ip<#1,#2>{\left\langle#1,#2\right\rangle}
\def\sp<#1>{\left\langle#1\right\rangle}

\newcommand{\eu}{\mathsf{eu}}
\newcommand{\Eu}{\mathsf{Eu}}
\def\ip<#1,#2>{\left\langle#1,#2\right\rangle}

\def\hp{\hphantom{x}}

\newcommand{\bplus}{\mbox{$\bigoplus$}}

\newcommand{\al}{{\alpha}}

\newcommand{\B}{{\mathfrak B}}
\newcommand{\A}{{\mathfrak A}}

\newcommand{\en}{{\enspace}}

\newcommand{\PP}{{\mathbb{P}}}

\def\oo{{\mathcal O}}
\def\K{{\mathcal K}}

\def\D{{\mathfrak{D}}}

\def\aa{{\scr A}}
\def\bb{{\scr B}}

\newcommand{\ta}{{\mathcal D}A }

\def\k{{\mathbb C}}

\def\Z{{\mathbb Z}}
\def\C{{\mathbb C}}

\begin{document}

\title{\en\en\large{Noncommutative del Pezzo surfaces\vskip 5pt
 and Calabi-Yau algebras }}

\author{Pavel Etingof}
\address{Department of Mathematics, Massachusetts Institute of Technology,
Cambridge, MA 02139, USA}
\email{etingof@math.mit.edu}

\author{Victor Ginzburg}\vskip 20mm
\address{Department of Mathematics, University of Chicago,
Chicago, IL 60637, USA}
\email{ginzburg@math.unchicago.edu}

\vskip 20mm
\begin{abstract}
The  hypersurface
 in $\k^3$ with an isolated quasi-homogeneous
 elliptic singularity of type
 $\wt E_r,r=6,7,8,$
  has a natural Poisson structure.
We show that the family of del Pezzo
surfaces of the corresponding type  $E_r$  provides
 a semiuniversal Poisson deformation of that  Poisson structure.

We also construct  a   deformation-quantization of the coordinate
ring of such a  del Pezzo
surface. To this end, we first deform the polynomial
algebra  $\k[x_1,x_2,x_3]$ to a  noncommutative
algebra   with  generators $x_1,x_2,x_3$
and  the following 3 relations labelled by cyclic parmutations  $(i,j,k)$ of
$(1,2,3)$:
$$
x_ix_j-t\cd x_jx_i=\Phi_k(x_k),\qquad
\Phi_k\in\C[x_k].
$$
This gives a family of Calabi-Yau  algebras $\A^{t_{\!}}(\Phi)$ parametrized by
 a  complex number $t\in\C^\times$ and a triple
$\Phi=(\Phi_1,\Phi_2,\Phi_3)$, of polynomials  of specifically chosen
degrees.

Our
quantization of the coordinate ring of a del Pezzo surface
is provided by noncommutative  algebras of the form
 $\A^{t_{\!}}(\Phi)/\llb\Psi\rrb$,
where  $\llb\Psi\rrb\sset
\A^{t_{\!}}(\Phi)$ stands for   the  ideal generated by a central element $\Psi$
which generates the center of the algebra $\A^{t_{\!}}(\Phi)$ if $\Phi$ is
generic enough.
\end{abstract}
\bigskip

\maketitle
\centerline{\sf Table of Contents}
\vskip -1mm

$\hspace{20mm}$ {\footnotesize \parbox[t]{115mm}{
\hp${}_{}$\!\hp1.{ $\;\,$} {\tt Introduction}\newline
\hp2.{ $\;\,$} {\tt Poisson deformations of a quasi-homogeneous surface
singularity}\newline
\hp3.{ $\;\,$} {\tt Main results}\newline
\hp4.{ $\;\,$} {\tt Three-dimensional Poisson structures}\newline
\hp5.{ $\;\,$} {\tt Poisson (co)homology}\newline
\hp6.{ $\;\,$} {\tt Classification results}\newline
\hp7.{ $\;\,$} {\tt Calabi-Yau deformations}\newline
\hp8.{ $\;\,$} {\tt From Poisson to Hochschild cohomology}\newline
\hp9.{ $\;\,$} {\tt{Appendix:\;computer computation}\qquad}
}}

\bigskip
{\large
\section{Introduction}\label{intro}}
\medskip

\subsection{Poisson structures on del Pezzo surfaces}\label{int2} We remind
the reader that a del Pezzo surface is a smooth projective surface
$S$ that is obtained by blowing up $\ell$
 sufficiently general points
in $\C{\mathbb P}^2$, where
$0\leq \ell\leq 8$, or $\C{\mathbb P}^1\times \C{\mathbb P}^1$. 
Let $S$ be such a   del Pezzo surface
with  canonical bundle $K_S$,
resp.  anti-canonical bundle $K\inv_S$.
A  regular section $\pi\in\Ga(S,K\inv_S)$
is a bivector  that gives $S$ a Poisson
structure (any bivector $\pi$ on a surface
automatically has a vanishing Schouten bracket:
$[\pi,\pi]=0$).
We say that a regular section $\pi\in\Ga(S,K\inv_S)$ is
{\em nondegenerate} provided the divisor of zeros of
$\pi$ is a reduced smooth curve.

In this paper we consider
 the most interesting case
where $\ell=6,7,$ or $8,$ and
where $\pi$ is assumed to be
a   nondegenerate  section.
Then, a simple application of the adjunction formula
shows that the zero locus of $\pi$ is
an elliptic curve $\E\sset S$. Furthermore,
$X:=S\sminus \E$ is an
affine surface equipped with an algebraic symplectic structure
provided by the (closed) 2-form $\pi\inv\in \Ga(S\sminus \E, K_S)$.

There are two Poisson algebras naturally  associated with
the data  $(S,\pi)$. The first algebra is
$\C[X]$, the coordinate ring of the
affine symplectic surface $X$.
The second algebra  is a graded algebra
\beq{RR}
\mR=\bplus_{n\geq0}\ \mR_n,\quad \mR_n:=\Ga\big(S,\ (K\inv_S)^{\o n}\big),
\eeq
the {\em homogeneous coordinate ring} 
associated with the anti-canonical bundle,
an ample line bundle on $S$.
One can use a construction of
Kaledin  to make $\mR$ a Poisson algebra as follows.

Choose a local nowhere vanishing section $\phi\in\Ga(U,\K_S\inv)$,
on a Zariski open subset $U\sset S$.
Let $L_h(\phi):=[i_\pi(dh), \phi]$ denote the Lie derivative
of the bivector $\phi$ with respect to  $i_\pi(dh)$, the Hamiltonian
vector field  associated with  a  regular function $h$ on $U$.
 Further, write $\{-,-\}_\pi$ for the Poisson bracket on
$U$ induced by the bivector $\pi$.
Then,  following \cite{Ka},
one defines
a Poisson bracket $\{-,-\}_\mR:\
\mR_n\times\mR_m\to\mR_{n+m},\ m,n\geq0,$ by the formula
$$\{f\phi^n,g\phi^m\}_\mR:=
\{f,g\}_\pi\cdot\phi^{n+m}+ \big(m g L_f(\phi)-
n f L_g(\phi)\big)\cdot \phi^{n+m-1},\qquad
\forall f,g\in\Ga(U,\oo_U).
$$
It is straightforward to verify  that the resulting bracket is
independent of the choice of a nowhere vanishing  section $\phi$, on $U$.

To relate the Poisson algebras $\C[X]$ and $\mR$, write $K$
for the total space of the canonical bundle $K_S$.
Thus  $K$ is a 3-dimensional variety equipped with a
natural $\C^\times$-action. By definition, one has
a graded algebra isomorphism $\mR=\C[K]:=\Ga(K, \oo_K)$,
such that the grading on $\C[K]$ comes from the
$\C^\times$-action. Further, there is a  diagram
\beq{diag1}
\xymatrix{
X=S\sminus \E\;\ar@{^{(}->}[rr]^<>(0.5){i:=\pi\inv}&&
K\ar@{->>}[rr]^<>(0.5){p}&& S,
}
\eeq
where the second map $p$ is the line bundle projection
and the first map
is a section of $p$ over $S\sminus \E$  provided by the symplectic form.

One can show  that the map $i=\pi\inv$,
in the diagram, is  a closed imbedding.
Moreover,
the corresponding  restriction
morphism $i^*:\ \mR=\C[K]\to \C[X]$
induces an algebra
isomorphism $\mR/(\pi-1)$ $\iso\C[X]$,
where $(\pi-1)$ denotes the ideal generated
by the element $\pi-1\in \mR_1\oplus\mR_0$.
The element  $\pi-1$ being {\em nonhomogeneous},
the grading on the algebra $\mR$
does {\em not} descend to a grading on
the quotient algebra. However, the ascending filtration
$F_{\leq m}\mR:=\bigoplus_{n\leq m} \mR_m$,
on $\mR$,
induces a well-defined  ascending filtration,
$F_\idot\C[X]$ that makes the coordinate ring $\C[X]$
a {\em filtered} algebra.

Let $\rees\C[X]:=\sum_{n\geq 0} F_n\C[X]\cdot t^n\sset \C[X]\o\C[t]$
be the {\em Rees algebra} of the filtered algebra
$\C[X]$. This is a graded algebra equipped
with a canonical graded algebra imbedding
$\C[t]\into \rees\C[X]$ such that one
has $\rees\C[X]/(t-1)\cong\C[X]$.
 Furthermore, the Poisson bracket on $\C[X]$
induces one on the Rees algebra.

We leave to the reader to prove the following simple result:

\begin{prop}\label{rees} There is a natural graded Poisson
algebra isomorphism  $\Xi:\ \rees\C[X]\iso \mR$ such that
\vskip2pt

\pb{The canonical  algebra imbedding
$\C[t]\into \rees\C[X]$ gets transported,
via $\Xi$, to the  graded algebra homomorphism
$\C[t]\into \mR$ induced by the assignment $t\mto\pi$.}
\vskip2pt

\pb{The isomorphism $\rees\C[X]/(t-1)\cong\C[X]$ gets transported,
via $\Xi$, to the algebra isomorphism $\C[K]/(\pi-1)\iso\C[X]$.\qed}
\end{prop}

The Proposition shows how the  Poisson algebras
$\C[X]$ and $\mR$ can be recovered from each  other.
Therefore, quantization (i.e., noncommutative
deformation) problems for these two algebras are
essentially equivalent.

In the rest of the
paper, we will concentrate on the problem
of quantizing the affine symplectic surface
$X=S\sminus\E$ by constructing  noncommutative
deformations of the Poisson algebra $\C[X]$, to be
viewed as coordinate rings of `noncommutative affine surfaces'.
Noncommutative
deformations of the  Poisson algebra $\mR$,
to be viewed  as homogeneous  rings of `noncommutative
{\em projective}
 surfaces',
may then be obtained by applying the Rees algebra construction
to the corresponding noncommutative deformations of  $\C[X]$.

\subsection{}\label{int1} The general theory of noncommutative
projective surfaces has been initiated in the late 80's
by Artin and Schelter \cite{AS}.
Many deep
results
were obtained later, in the papers \cite{ATV}, \cite{AV}, \cite{BSV},
\cite{C1}, \cite{Le},
and
\cite{St1}.

The general philosophy of  noncommutative
 surfaces, either projective or affine, was outlined by M. Artin in
\cite{A}.
According to that philosophy in the
{\em affine} case, one tries to construct
a noncommutative  algebra $B$ that plays the role of
`coordinate ring' of an  (affine) noncommutative surface $X$.
It turns out that, in typical examples,
 the algebra $B$ often appears   in the form
$B=A/\llb\Psi\rrb$. Here, $A$ is an auxiliary associative algebra
which is somehow more accessible than $B$, and $\llb\Psi\rrb$
denotes a two-sided ideal in $A$ generated by
 a normal (often central)
element $\Psi\in A$. It has been
remarked by M. Artin  \cite{A} that there should be some  more conceptual
 {\em a priori} explanation  of the appearance   of the
algebra $A$ and  of the element
$\Psi$.

The aim of the present paper is to propose such an explanation. Our
approach is based on the concept of Calabi-Yau  (CY) algebra, introduced recently,
cf. eg. \cite{Bo}, \cite{Gi}, and Definition \ref{vdb_def}
below. This approach is consistent with the
point of view of string theory where 3-{\em dimensional} CY varieties
are considered  to be more fundamental  than 2-dimensional surfaces.
Thus, a  2-dimensional surface should be viewed as a hypersurface
in an ambient CY 3-fold which, in the affine case, is typically taken to be
 $\k^3$ and, in the projective case,
is taken to be the total space of the canonical bundle
of the surface.

The best way to understand what kind of
noncommutative algebraic structers should be analogous to the
structures of CY geometry is to consider a `quasi-classical approximation' first.
A noncommutative  CY algebra of dimension 3 reduces,
quasi-classically,  to the coordinate ring $\k[M]$ of an
  affine 3-dimensional variety $M$. Such a variety  comes equipped with
  an
algebraic volume
form $\vol\in\Om^3(M),$ that  keeps track of the CY structure, and
with a Poisson bracket,
that `remembers' about the noncommutative deformation, up to first order.
A key point is  that these two pieces of data must be related. Specifically,
it was explained by Dolgushev \cite{Do} that the correct
quasi-classical analogue of the CY condition is the
requirement that the  Poisson bracket on $M$ be {\em unimodular},
that is, such that any Hamiltonian vector field on $M$ preserves
the volume form $\vol$, i.e. has the vanishing divergence.

It is easy to show that any  unimodular  Poisson bracket on
a 3-fold with trivial first de Rham cohomology is determined
by a single regular function $\phi\in\k[M]$, see \S\ref{3}.
The function $\phi$ is unique up
to a constant summand and it is automatically central with respect to the corresponding
Poisson bracket. Furthermore, this function generates, generically, the
whole Poisson center.

We turn now to noncommutative surfaces inside our noncommutative CY variety.
Quasi-cla- ssically, giving such a surface amounts
to giving a {\em Poisson} hypersurface $X\sset M$.
 For  $M=\k^3,$ for instance,  that means, in the generic case,
that the equation of the  hypersurface $X$ must be given
by a function contained in the Poisson center of $\C[M]$.
In the situation where the  Poisson center reduces to
$\C[\phi]$ we conclude that
our function is a polynomial
in $\phi$. Hence, the only  hypersurfaces  which may arise
in the process of  quasiclassical degeneneration of
a noncommutative story
are, essentially,
the level sets of $\phi$. By redefining $\phi$,
one may assume without loss of generality that
the surface is  the zero set of $\phi$, so
the corresponding coordinate
ring is $\C[X]=\C[M]/(\phi)$.

The discussion above suggests that $\C[M]$, the
coordinate ring of the CY  3-fold, gets
deformed via a quantization to a noncommutative CY algebra $A$ in such
a way that the  function $\phi$ gets deformed to
a {\em central} (more generally, normal) element
$\Psi\in A$.  Therefore, the coordinate ring of the
corresponding surface  gets deformed
to a noncommutative algebra of the form $B=A/\llb\Psi\rrb.$

This provides a reason for the appearance of
 the objects $A$ and $\Psi$ we were looking for.

\subsection{}
In this paper, we study hypersurfaces in the CY variety
 $M=\k^3$,  equipped with the standard
volume form $\dr x\wedge\dr y\wedge\dr z$.
Thus, we have $\C[M]=\k[x,y,z]$.
As we have mentioned earlier,
associated with any  $\phi\in\C[x,y,z]$, there is
a  Poisson structure on $M$. Specifically,
the Poisson
brackets of
coordinate functions are given by the following explicit formulas
\beq{formula}
\{x,y\}=\frac{\partial \phi}{\partial z},\quad
\{y,z\}=\frac{\partial \phi}{\partial x},\quad
\{z,x\}=\frac{\partial \phi}{\partial y}.
\eeq

It is immediate to verify that   $\phi$ is a central element
with respect to the above bracket. Therefore,
$\k[x,y,z]/(\phi)$,  a quotient by the principal ideal generated by $\phi$,
 inherits  the structure of a Poisson
algebra.

\begin{defn}\label{aabb} We write  $\aa_\phi:=\k[x,y,z]$ for the Poisson algebra
with bracket \eqref{formula}, and let $\bb_\phi:=\aa_\phi/(\phi)$
be the quotient Poisson algebra with induced bracket.
\end{defn}

It is interesting to take $\phi$  a (quasi-) homogeneous
polynomial  with an isolated singularity at the origin.
In the special case where $\deg\phi\leq \deg x+\deg y+\deg z$,
the equation $\phi=0$ defines  a Poisson surface
with either  simple Kleinian, or elliptic
singularity.

We study both commutative and noncommutative
 deformations
of the corresponding Poisson algebra $\bb_\phi$.
We show that
all Poisson algebra deformations  are essentially obtained
by deforming the polynomial $\phi$, see Theorem \ref{delpezzoform}.
In the elliptic case, for instance,
any such deformation  gives the coordinate ring of an affine surface
obtained by removing an elliptic
curve from an appropriate projective del Pezzo surface.

Our approach to  noncommutative
deformations  of  elliptic singularities is motivated by the ideology explained
in \S\ref{int1}. Specifically, we {\em  simultaneously} deform both the
corresponding surface $\phi=0$
and the ambient CY variety $\k^3$.
This way, we construct a flat family of
noncommutative CY  algebras $\A^{t^{\!}}(\Phi)$ of dimension 3,
which provide
a deformation of the  Poisson  algebra $\aa_\phi$,
and a family of central elements $\Psi\in \A^{t^{\!}}(\Phi)$.
The  noncommutative algebras of the form $\A^{t^{\!}}(\Phi)/\llb\Psi\rrb$
thus provide a flat deformation
 of
the  Poisson  algebra $\bb_\phi$.
In analogy with the Poisson case, these noncommutative algebras may be thought of as
`coordinate rings' of noncommutative  del Pezzo surfaces.

There were a few other approaches to the problem of
 quantization of   del Pezzo
surfaces in the literature. One of them was 
proposed by M. Van den Bergh, in the paper \cite{VB3}, 
which gives a construction of the category of coherent sheaves
on a `would be' noncommutative (projective)  del Pezzo surface.
The connection between this approach and our approach is given by Chapter 12 of \cite{VB3}. 
Namely, it is shown there that if one blows up 6 points in a quantum plane and then takes the affine part 
(the complement of the elliptic curve), then the coordinate ring is of the form $A/(n)$, 
where $A$ is a filtered deformation of an AS-regular algebra and $n$ is a normalizing element. 
We expect that this ring is the $E_6$-deformation considered in this paper,
and that a similar approach works for $E_7$ and $E_8$. 

 A different construction which is explicit but works only for
a  very special class of  {\em degenerate}  noncommutative  del Pezzo
 surfaces, was proposed in \cite{EOR}.

Our present approach works in the general case, and is both quite simple and   explicit.
As a first step, we introduce a family of associative algebras
$\A^{t_{\!}}(\Phi)$ to be the algebras with 3 generators, $x,y,z$,
subject to 3 defining relations of the following  form
\beq{for1}
[x,y]_t=\frac{\partial \Phi}{\partial z},\quad
[y,z]_t=\frac{\partial \Phi}{\partial x},\quad
[z,x]_t=\frac{\partial \Phi}{\partial y},
\eeq
In this formula,
$\Phi$ runs over  a certain explicitly defined
family of noncommutative {\em cyclic potentials},
 $t$ is a complex parameter, and we have used the notation $[u,v]_t:=uv-t\cdot vu$.

\begin{rem}
It is interesting to note  that relations in \eqref{for1}
look very similar to the formulas for the Poisson bracket
\eqref{formula},
at least formally. The analogy  goes much further
since the actual formula for $\Phi$, see \eqref{xyz}-\eqref{XYZ},
is quite similar to the formula for the polynomial
$\phi\in\k[x,y,z]$ that gives the equation of an
affine del Pezzo surface, see \eqref{pol}-\eqref{phiPQR}.
\end{rem}

Next, we prove one of our main results, see \S\S\ref{AABB}-\ref{NCDP},
saying that
$\A^{t_{\!}}(\Phi)$ is a Calabi-Yau algebra of dimension 3
and that,   for sufficiently general parameters,
 the center of $\A^{t_{\!}}(\Phi)$ has the form $\k[\Psi]$,
a free polynomial algebra  generated by an element $\Psi$ uniquely
determined up to a constant summand.
We  show further that the family of
noncommutative algebras of the form $\B^{t_{\!}}(\Phi,\Psi):=\A^{t_{\!}}(\Phi)/\llb\Psi\rrb$
provides the required quantization of del Pezzo surfaces.
It is also quite remarkable that, in a sense, {\em any}
flat infinitesimal deformation of the Poisson algebra
$\bb_\phi$ can be obtained by the above construction,
cf. Theorem \ref{two}.

In section \ref{NCES}, we discuss the special case of homogeneous
potentials. In this case, the algebras  $\A^{t_{\!}}(\Phi)$ and
$\B^{t_{\!}}(\Phi,\Psi)$ have natural gradings.
The graded algebra $\A^{t_{\!}}(\Phi)$
is nothing but an Artin-Schelter regular algebra
of dimension 3. These algebras, also known as
{\em Sklyanin algebras}, have been intensively studied
in the literature, see \cite{AS}, \cite{ATV}, \cite{AV} and references
therein. In particular, they were classified in D. Stephenson's Ph.D. thesis, \cite{St2}
(see also \cite{St3}). 
The best understood case is that of singularities of type $E_6$, resp. $E_7$,
corresponding to quadratic, resp. cubic, Sklyanin algebras.
The $E_8$-case hasn't been studied so well, cf. however ~\cite{St1}.

The graded algebra $\B^{t_{\!}}(\Phi,\Psi)$
may be thought of as the homogeneous coordinate ring of a
noncommutative  elliptic singularity.
There seems to be an interesting and largely unexplored theory of
graded {\em matrix factorizations} for noncommutative
elliptic singularities.  In section \ref{MATR}, we introduce a few
basic results, cf. also \cite{KST}, and formulate Conjecture
\ref{matr_conj}.

In the general case of  an arbitrary,
not necessarily homogeneous, potential $\Phi$
the algebra  $\A^{t_{\!}}(\Phi)$ comes equipped
with a natural ascending filtration and
one may form the  corresponding {\em Rees algebra}.
This way, one obtains
 a class of graded algebras that has been considered
earlier,
especially in  type $E_6$, see  \cite{BSV} and \cite{C1}, \cite{C2}.
Nonetheless,  an explicit expression
for the central  element $\Psi\in \A^{t_{\!}}(\Phi)$,
or the corresponding homogeneous central element
of the Rees algebra,  is quite
complicated even in type $E_6$, see \S\ref{magma_sec} and \cite{R}.

\begin{rem}
It would be  interesting to establish a
connection between our approach to noncommutative
del Pezzo surfaces and the results of Chan-Kulkarni \cite{CK}.
\end{rem}

\subsection{Definition of Calabi-Yau algebras}
 We will work with unital
associative, not necessarily commutative, $\C$-algebras,
to be referrred to as `algebras'. We write $\o=\o_\k,\,\dim=\dim_\k,$
etc.

\begin{defn}[\cite{Gi}]\label{vdb_def} An
algebra  $A$ is said to be a {\em Calabi-Yau algebra
of dimension} $d\geq 1$, provided it has  finite Hochschild dimension, and
there are  $A$-bimodule isomorphisms
\beq{cy_def}
\Ext_{\bimod{A}}^k(A,A\otimes A)\cong
\begin{cases} A &\operatorname{if}\enspace k=d\\
0&\operatorname{if}\enspace k\neq d.
\end{cases}
\eeq
 \end{defn}

The image of  $1_A\in A$
under such an isomorphism gives a {\em central} element in
 $\Ext_{\bimod{A}}^d(A,A\otimes A),$
called   {\em noncommutative volume} on $A$.

\begin{examp} Let $X$ be a smooth connected affine complex algebraic
variety of dimension $d$. A noncommutative volume for the algebra
$A=\k[X]$, the coordintate ring of $X$, is the same thing as
a nowhere vanishing section
of the line bundle $\La^d T_X={{\scr E}^{\!}xt}^d_{\oo_{X\times
X}}(\oo_X,\oo_{X\times X}).$ Thus,
 $A$ is a Calabi-Yau algebra  if and only if  $X$
is a Calabi-Yau variety.
\hfill$\lozenge$\end{examp}

\begin{rem} Following Van den Bergh \cite{VB1}, it may be natural
to consider a wider class of
{\em twisted} CY algebras which satisfy a
weaker version of  \eqref{vdb_def} requiring that the group
$\Ext_{\bimod{A}}^k(A,A\otimes A)$ be zero for $k\neq d$ and,
 for $k=d,$ this Ext-group be an {\em arbitrary invertible}
$A$-bimodule $U$, not necessarily $U=A$.
Twisted CY algebras correspond geometrically to
 arbitrary Gorenstein varieties whose dualizing
sheaf is a not  necessarily trivial line bundle.

One should be able to develop an analogue of the theory of CY algebras
in this more general framework. In such a theory, the role
of  $\dr\Phi$, an {\em exact} noncommutative cyclic 1-form associated with
a cyclic potential $\Phi,$ cf. \S\ref{defs} and \cite{Gi} \S 3.5, is
expected to be played by a suitable noncommutative cyclic 1-form
with coefficients in $U\inv$, an inverse $A$-bimodule.

In the special case of graded algebras, any invertible graded
$A$-bimodule $U$ must be a rank 1 free left $A$-module. The right
$A$-action on $U$ is
then
given, in terms of a left $A$-module isomorphism $U\cong A,$  by the formula
$ua=u\cdot\sigma(a),\forall u\in U,a \in A,$
where $\sigma$ is an algebra automorphism of $A$.
In the framework of Sklyanin algebras, this
 has the effect
that the {\em central} element $\Psi$ of the CY algebra gets replaced by a {\em normal}
element in a twisted CY
algebra, cf. \cite{ATV}.
\end{rem}

\subsection{Acknowledgements}{\footnotesize
We are  grateful to Mike Artin and  Eric Rains for useful discussions.
We also thank Eric Rains for
allowing us to reproduce his computer
computations in the appendix to this paper.
We are grateful to the referee for comments and references. 
The work of  E.G. was  partially supported by the NSF grant DMS-0504847.
The work of V.G. was  partially supported by the NSF grant DMS-0601050.}

\medskip
{\large
\section{
Poisson deformations of a quasi-homogeneous surface
singularity}\label{intro_sub}}
\medskip
\subsection{Deformations and  cohomology}\label{basic}
Deformations of an algebraic  object $A$
are often
controlled by the vector space
$H^2(A)$, the second cohomology group for an appropriate cohomology theory.
That means, in particular, that
associated with such a deformation, i.e. with a family of
objects $\{A_s, s\in S\}$ parametrized by a scheme $S,$
one has  a canonical classifying,
Kodaira-Spencer type, linear map
\beq{KS}
\KS_s:\
T_sS\to H^2(A_s),\qquad s\in S,
\eeq
where $T_sS$ stands for the Zariski tangent space to the scheme $S$ at a
point $s$.
A tangent vector
$v\in T_sS$ determines a 1-parameter infinitesimal first order
deformation of the object
 $A_s$. The image of $v$ under the classifying map $\KS_s$ is called the Kodaira-Spencer
class of that deformation.

\begin{defn} A  family  $\{A_s, s\in S\}$,
parametrized by a {\em smooth} scheme $S,$ is  said to be a (smooth) {\em semiuniversal}
\footnote{The term ``semiuniversal deformation'' is often used for deformations parametrized by arbitrary 
(not necessarily smooth) formal schemes. In this paper, we will consider only smooth semiuniversal deformations, and 
for this reason will not explicitly mention that they are smooth.}
deformation provided the classifying map is a vector space isomorphism for any $s\in S$.
\end{defn}

Obstructions to deformations of an object $A$ are often
controlled by
$H^3(A)$, the third cohomology group.
 A standard result of deformation theory insures the existence
of  a formal semiuniversal   deformation of
 $A$ with base   $S=H^2(A)$
provided  one has: \en\textsf{(1)} $\dim H^2(A)<\infty$ and, moreover,
\en\textsf{(2)} $H^3(A)=0$. However, a formal semiuniversal deformation of $A$ sometimes exists 
even if $H^3(A)\ne 0$. If the semiuniversal deformation exists, one says that the deformations of $A$ are unobstructed.

Given an associative, resp. commutative associative or Poisson, algebra $ A $, one can  define its
Hochschild cohomology $\HH^\hdot(A):=\Ext_{\bimod{A}}^\hd(A,A)$ (Gerstenhaber), resp. Harrison
cohomology,
$\harr^\hd(A)$ (cf. \cite{Lo} and references therein),
 or  Poisson cohomology $\PH^\hdot( A )$
(cf.  \cite[Appendix]{GK} and \S\ref{hp_subsec} below).
By definition, in degree zero for an associative algebra $A$ one has
$\HH^0(A)=Z(A),$  the center of  $A$.
Similarly, for a Poisson algebra with
 Poisson bracket $\{-,-\}: \ A\times A\to A$, we have
 $\PH^0(A)=\zz(A):=\{z\in A\,\mid\, \{z,a\}=0,\en\forall a\in A\}$,
is the {\em Poisson
 center} of  $A$.

Also, for the corresponding degree zero Hochschild, resp.  Poisson,
{\em homology}, one has $\HH_0(A)=A\nat:=A/[A,A],$
the commutator quotient space, resp.
$\PH_0(A)=A/\{A,A\}.$

 Flat deformations of an associative, resp. commutative associative
 or Poisson, algebra $ A $
are controlled by the second
Hochschild cohomology group
 $\HH^2(A),$ resp.
$\harr^2(A)$ or   $\PH^2( A ),$
cf. \cite{GK}.
Thus, one may consider {\em flat}  deformations of such
an algebra $A$. Observe that,
 a  flat family of
 Poisson algebras is in particular
 a  flat family of commutative algebras.
This corresponds, in terms of cohomology,
to the existence for any  Poisson algebra $A$ of a canonical
linear map
$\textsf{can}:\
\PH^2(A)\to\harr^2(A)$.

Now, let $A$ be a Calabi-Yau algebra  of dimension $d$ in the sense
of Definition \ref{vdb_def}.
According to \cite{VB2}, a choice of noncommutative  volume  for
 $A$ induces
a Poincar\'e duality type isomorphism
\beq{vdb}\HH_\idot(A) \iso \HH^{d-\hdot}(A).
\eeq

 Following \cite{CBEG}, we introduce
a BV operator $\De: \HH^\hdot(A)\to \HH^{\hdot-1}(A)$,
 obtained by transporting the Connes differential $\BB$, on
Hochschild homology, to Hochschild {\em co}homology
via the duality isomorphism \eqref{vdb}.

One may consider first order deformations of the CY algebra $A$ within
the class of Calabi-Yau  algebras.
The
Kodaira-Spencer classes of all such deformations form a  vector
{\em sub}space in  $\HH^2(A),$ that  turns out to be equal to
$$ \Ker[\De: \HH^2(A)\to\HH^1(A)].
$$

In the special case of Calabi-Yau  algebras of dimension $d=3$,
there is a
chain of maps
\beq{kappadeform}
\kappa:\
\xymatrix{
A\nat\ar@{->>}[rr]^<>(0.5){\eqref{cy_def}}_<>(0.5){\sim}&&\HH^3(A)\
\ar[r]^<>(0.5){\Delta}
&
\en\Ker[\De: \HH^2(A)\to\HH^1(A)],
}
\eeq
where we have used that $\Image(\De)\sset\Ker(\De),$ since $\De^2=0$.

Let $A=\A(\Phi)$ be a Calabi-Yau  algebra of dimension 3 defined by a potential
$\Phi$, see \S\ref{defs}.
An arbitrary infinitesimal variation
$\Phi\rightsquigarrow \Phi+\eps\Phi'$ (where $\eps^2=0$), of the potential,
yields an infinitesimal deformation
of
 $A$. We show in \S\ref{def_sub} below, cf. also \cite{BT},  that such a deformation is
 automatically
flat; moreover, it is a deformation within
the class of Calabi-Yau  algebras.
Let $\Phi'\nat\in A\nat$ denote the class of $\Phi'$
in the commutator quotient.
Then, it is not difficult to prove the following proposition, whose proof is left to the reader. 

\begin{prop}\label{KS_prop}
The Kodaira-Spencer class in $\Ker[\De: \HH^2(A)\to\HH^1(A)]$
of the deformation  $\A(\Phi)\rightsquigarrow
\A(\Phi+\eps\Phi')$ is equal to
$\kappa(\Phi'\nat),$ the image of $\Phi'\nat$ under the
composite map \eqref{kappadeform}.
\qed
\end{prop}

\subsection{Quasi-homogeneous surface singularities.}\label{coh}
Let the multiplicative group $\C^\times$ act on $\k^3$
with positive integral weights  $a\leq b\leq c.$
This makes  the coordinate ring
$\k[x,y,z]$, of $\C^3$, a nonnegatively graded  algebra
with homogeneous generators of degrees
$\deg x=a,\,\deg y=b,\,\deg z=c.$
Thus,
 $\phi\in\k[x,y,z] $ is a
(weighted-, equivalently, quasi-) homogeneous  polynomial of weight $\deg \phi=\d$
if and only if one has $\eu(\phi)=\d\cdot \phi,$ where
\beq{eu}
\eu:=ax\frac{\partial}{\partial x}+
by\frac{\partial}{\partial y}+
cz\frac{\partial}{\partial z},
\eeq
denotes the {\em Euler vector field} that generates the $\C^\times$-action.

Associated with any polynomial $\phi\in \k[x,y,z]$ with an {\em isolated}
singularity, is its  {\em Jacobi ring}
 $\J(\phi):= \k[x,y,z]/(\frac{\partial   \phi}{\partial x},
\frac{\partial   \phi}{\partial y},
\frac{\partial   \phi}{\partial z}).$
If  $\phi$ is (weighted-) homogeneous of weight
 $\d,$ then $0\in\C^3$ is the only singular point.
 Furthermore, the  Jacobi ring acquires a natural
grading $\dis\J(\phi)=\oplus_{m\geq 0}\,\J^{(m)}(\phi).$
 For the corresponding Hilbert-Poincar\'e polynomial,
one easily finds the formula, cf. \S\ref{spectr},
\beq{R}
\sum\nolimits_{m\geq 0} u^m\cdot\dim\J^{(m)}(\phi)=
\ \frac{(u^{\d-a}-1)(u^{\d-b}-1)(u^{\d-c}-1)}{(u^a-1)(u^b-1)(u^c-1)}.
\eeq

Set  $\mm_\phi:=\phi\inv(0)\sset\k^3$.
Specializing the RHS of \eqref{R} at $u=1$, we get a formula
\beq{mu}
\dim \J(\phi)=\mu:=\frac{(\d-a)(\d-b)(\d-c)}{abc},
\eeq
for  the  {\em Milnor number} of the isolated singularity (at the origin) of
 the hypersurface  $\mm_\phi$.

 Let  $a\leq b\leq c<\d$ be an arbitrary quadruple of
positive integers
such that
${\mathsf{gcd}}(a,b,c,\d)=1.$ According to Kyoji Saito  \cite[Theorem
3]{Sa},
one has the following result.

\begin{thm}[Saito]\label{saito} Assume that  the rational function
 associated with the  quadruple $(a,b,c;\d)$
by the formula on the right of \eqref{R}
 is a polynomial (i.e. has no poles).

Then,
 the
surface $\mm_\phi$  has an isolated
singularity at the origin, for any general
enough homogeneous polynomial $\phi\in\k[x,y,z],$ of
degree
$d$.\qed
\end{thm}

\subsection{Simple Kleinian and elliptic singularities}\label{simple}
Let $\PP^{a,b,c}=(\k^3\sminus\{0\})/\C^\times$ denote the  weighted projective
plane corresponding  to the $\C^\times$-action with weights $(a,b,c)$,
where $\text{gcd}(a,b,c)=1.$
Restricting the   projection
$\k^3\sminus\{0\}\to\PP^{a,b,c}$ to the punctured  hypersurface, one
obtains a map $\mm_\phi\sminus\{0\}\onto \PP(\mm_\phi)\sset \PP^{a,b,c}$.
This way   $\mm_\phi\sminus\{0\}$ becomes
a principal $\C^\times$-bundle over $\PP(\mm_\phi)$, a projective
curve.
The type  of the hypersurface $\mm_\phi$
is closely related to the integer
\beq{varpi}
\pp:=\d-a-b-c.
\eeq

There is a complete list
of all   hypersurfaces with $\pp=-1,0,1,$ see \cite{Sa}.
According to K. Saito, for any such   hypersurface,
one has $\mm_\phi\sminus\{0\}\cong H_\pp/\Gamma$.
Here, $H_\pp$ is the total space of the
$\C^\times$-bundle associated with
the   canonical line bundle on  a curve $C_\pp$, and
$\Gamma$ is a discrete group of bundle automorphisms.
Depending on whether  $\pp=-1,0,$ or $+1,$
the curve $C_\pp$ is either
the projective line $\PP^1(\C)$,
or
the affine line $\C$, or  the upper half plane,
 respectively.
Moreover, in each case, the  group $\Gamma$ is a discrete
subgroup of the group of motions of  $C_\pp$,
viewed as a Riemann surface with the natural metric,
and the $\Gamma$-action on $H_\pp$ is induced by
the natural  $\Gamma$-action on $C_\pp$.

In  the case $\pp=-1,$ the surface $\mm_\phi$ has a
simple $A,D,E$ (Kleinian)
singularity, while  the case
$\pp=0$ corresponds to simple elliptic singularities
$\widetilde E_6,\,\widetilde E_7,\,\widetilde E_8$
(for a  reducible curve, all components
must be rational). Specifically, one has
the following classical result, cf.
e.g. \cite{B}, and \S\ref{class} below.

\begin{prop}\label{alggeo}
Let the variables $x,y,z$ have degrees $0<a\leq b\leq c$, such that
${\rm gcd}(a,b,c)=1$.

\vi Let
 $\phi\in\C[x,y,z]$ be an  irreducible  homogeneous polynomial
 of degree $\deg\phi\le a+b+c$. Then, the projective curve $\phi(x,y,z)=0$ is either rational or elliptic.
\vskip 2pt

\vii Let $\d\le a+b+c$ be such that, for a general  homogeneous polynomial
$\phi$ of degree $\d$, the  projective curve  $\phi(x,y,z)=0$ is elliptic.
Then,  $\d=a+b+c,$ and we have
\begin{itemize}

\item One of the
following   holds:

\beq{table}
{\renewcommand{\arraystretch}{1.2}
\begin{tabular}{ll|lll|c|ccc|c}
  &                         & a& b& c& d\en
  &$
\dis p:=\fra{d}{a_{}}$&$\dis q:=\fra{d}{b_{}}$&$\dis r_{_{}}:=\fra{d}{c_{}}$
&{$\mu$}\\
\hline
$1.$& $E_6\en\oper{case}$\,:\ & $1$& $1$& $1$&$ 3$\en &$3$&$3$&$3$&$8$\\
$2.$& $E_7\en\oper{case}$\,:\ & $1$& $1$& $2$&$ 4$\en &$4$&$4$&$2$&$9$\\
$3.$& $E_8\en\oper{case}$\,:\ & $1$& $2$& $3$&$ 6$\en &$6$&$3$&$2$&$10$
\end{tabular}
}
\eeq
moreover, the integers $\dis(p-1,
q-1, r-1)$ give the
lengths of 3 legs of the corresponding {\em extended}
Dynkin diagram of type $\wt E_6,\wt E_7,$ or $\wt E_8$.

\item The homogeneous equation of the corresponding elliptic curve
can be brought to the canonical form
\beq{pqr}
\phi^\tau(x,y,z)=\frac{x^{p}}{p}+\frac{y^{q}}{q}+ \frac{z^{r}}{r} + \tau\cd xyz=0,\quad\oper{where}\quad
 \tau\in\C^\times.
\eeq
\end{itemize}
\end{prop}

We note that, in the setting of \eqref{table} one has
$\ \dis
\fra{1}{p}+\fra{1}{q}+\fra{1}{r}=\fra{a}{d}+\fra{b}{d}+\fra{c}{d}
=\fra{a+b+c}{d}=1.
$

\begin{rem}
The case
 $\pp=1$ turns out to be closely related to
14 exceptional  singularities (Dolgachev singularities)  arising
in degenerations of $K3$ surfaces.
\end{rem}

\subsection{}
Let $\C^\times$ act on  $\k^3$ with  weights
$0< a\leq b\leq c$, where ${\mathsf{gcd}}(a,b,c)=1.$
Associated with $\phi\in\k[x,y,z]$ we have
 the Poisson  algebra $\bb_\phi$, see Definition \ref{aabb}.

The following theorem will be proved
in Subsection \ref{Poihom}
 using some results of  Pichereau, \cite{P}, 
explained in \S\ref{pcoh}.

\begin{thm}\label{compute} For  a (quasi-) homogeneous
polynomial $\phi$ with an isolated singularity, we have

\vi The   Hochschild cohomology of $\bb_\phi$ is as follows
$$\HH^\hdot(\bb_\phi)\cong\X^\hdot\bb_\phi
\bigoplus u^2\cd \k^\hdot[u]\o \J(\phi),
\qquad\deg u=1.$$

\vii  The Poisson
cohomology  of $\bb_\phi$
is as follows
$$
\PH^0(\bb_\phi)=\k,\quad
\PH^1(\bb_\phi)=\J^{(\pp)}(\phi),\quad
\PH^2(\bb_\phi)=\J^{(\pp)}(\phi)\oplus \J(\phi),
\quad\PH^{k}(\bb_\phi)=\J(\phi),\en k\geq 3.
$$
\end{thm}

Here, in part (i), $\X^\hdot\bb_\phi$ denotes the algebra of poly-derivations
of the algebra $\bb_\phi,$  cf. \S\ref{not}, and in part
(ii) we use the notation \eqref{varpi}.

Theorem \ref{compute}(ii) shows that the group $PH^3(\bb_\phi)$
does not vanish.
Nonetheless, 
there is an explicit 
Poisson deformation of the Poisson algebra $\bb_\phi$ such that the
tangent space to the base of that deformation
is identified with  $\PH^2(\bb_\phi)=\J^{(\pp)}(\phi)\oplus \J(\phi)$.
Specifically, the space $\J(\phi)$, the second 
direct summand, parametrizes
 deformations of the Poisson algebra
 $\bb_\phi$ obtained by
deformations of the polynomial
$\phi$. Any nontrivial deformation of this kind
gives a nontrivial deformation of  $\bb_\phi$,
viewed as a commutative algebra (with the Poisson structure disregarded),
cf. relation to Harrison cohomology below.

On the other hand, the space $\J^{(\pp)}(\phi)$,  the first direct summand
in the decomposition
$\PH^2(\bb_\phi)=\J^{(\pp)}(\phi)\oplus \J(\phi),$
parametrizes deformations which change the
Poisson structure on $\bb_\phi$ while keeping
the commutative algebra structure unaffected.
To see this, we use results of Pichereau \cite{P},
see also formula \eqref{piphi} in the paper.
According to {\em loc cit},
elements of the direct summand  $\J^{(\pp)}(\phi)\sset
\PH^2(\bb_\phi)$ may be represented
by bivectors of the form $f\cdot\pi,$
where $f\in \bb_\phi$ is a homogeneous
element of degree $\pp$ and
$\pi$ is the Poisson bivector that gives
the Poisson bracket \eqref{formula}, on $\bb_\phi$.
We will see in the course of the
proof of Theorem \ref{compute} that the family
of bivectors of the form $\pi+f\cdot\pi,\ f\in\J^{(\pp)}(\phi),$
yields the required family of nontrivial deformations
of the Poisson structure on $\bb_\phi$, parametrized
by the vector space  $\J^{(\pp)}(\phi)$.

The direct sum
decomposition $\HH^2(\bb_\phi)=\X^2\bb_\phi\oplus u\cdot\J(\phi),$
in Theorem \ref{compute}(i),
corresponds to the Hodge decomposition of Hochschild
cohomology, cf. \cite{Lo}, \S 4.5.
The second direct summand is equal to  $\harr^2(\bb_\phi),$
the second  Harrison cohomology group of 
the algebra $\bb_\phi$.
By general deformation theory, the latter group
is the base
 of the semiuniversal
unfolding of the quasi-homogeneous isolated singularity
$\phi=0$. Thus,
the canonical morphism $\textsf{can}:\
\PH^2(\bb_\phi)\to  \harr^2(\bb_\phi)$,
that send a Poisson deformation to the
corresponding deformation of the underlying commutative algebra,
may be identified with the second projection
$\J^{(\pp)}(\phi)\oplus \J(\phi)
\to \J(\phi)$. This agrees with the discussion of the preceding paragraph:
the direct summand that corresponds 
to Poisson deformations of $\bb_\phi$ induced by deformations of
the polynomial $\phi$ projects isomorphically onto the group
$\harr^2(\bb_\phi)$. On the other hand, the direct summand
$\J^{(\pp)}(\phi)$, that corresponds 
to  deformations of the Poisson structure which do not change
the commutative algebra structure, 
projects to zero.

Note that if $\pp=-1$, the
case of Kleinian singularity, Theorem \ref{compute} yields
$\PH^1(\bb_\phi)=0$ and
$\PH^2(\bb_\phi)=\harr^2(\bb_\phi)=\J(\phi)$.
It is easy to see that,
in this case, the map $\textsf{can}$ reduces to the identity.

\subsection{Poisson deformations of elliptic singularities} \label{ellip}
Given a triple $(p,q,r)$ of positive integers, introduce a triple
of  polynomials
\begin{align}
P&=\fra{1}{p}\cd x^p+ \alpha_1\cd x^{p-1}+\ldots+\alpha_{p-1}\cd x\in \k[x],\nonumber\\
Q&=\fra{1}{q}\cd y^q+\beta_1\cd y^{q-1}+\ldots+\beta_{q-1}\cd y\in \k[y],\label{pol}\\
R&=\fra{1}{r}\cd z^r+\gamma_1\cd z^{r-1}+\ldots+\gamma_{r-1}\cd z\in \k[z].\nonumber
\end{align}

 Further, we let
\beq{phiPQR}
\phi^{\tau,\nu}_\pqr:=\tau\cd xyz+P(x)+Q(y)+R(z)+\nu\in \k[x,y,z],\qquad\
\tau\in\C^\times,\,\nu\in\C.
\eeq
The family of  polynomials $\phi^{\tau,\nu}_\pqr$ depends on
$(p-1)+(q-1)+(r-1)+2=p+q+r-1$
complex parameters $\alpha_i,\beta_j,\gamma_k,\tau,\nu$.
If all the parameters, except for the parameter $\tau$, vanish,
this family
specializes to a homogeneous polynomial
$\phi^\tau=\phi^{\tau,0}_{0,0,0}$
of the form
\eqref{pqr}.

Recall that, for any polynomial $\phi\in\C[x,y,z]$,
the equation $\phi(x,y,z)=0$ defines an
 affine Poisson  surface  in $\k^3$,
with  coordinate ring $\bb_{\phi}$.

\begin{thm}\label{delpezzoform} Let $(a,b,c)$ and $(p,q,r)$ be the integers associated to
 one of the 3 cases $E_\ell, \ell=6,7,8,$ of
table \eqref{table}, and let $\phi^\tau$ be the corresponding
polynomial \eqref{pqr}. Then,
\vskip 2pt

\vi  For the Milnor number $\dim\J(\phi^\tau)=\mu$, we have
\beq{dimJ}
\mu=\fra{(a+b)(a+c)(b+c)}{abc}=p+q+r-1.
\eeq

\vii The equations $\phi^{\tau,\nu}_\pqr(x,y,z)=0$ give a flat
$\mu$-parameter family of  affine del Pezzo
surfaces of the corresponding type $E_\ell,\,\ell=6,7,8.$
\vskip 3pt

\viii The family of
Poisson algebras $\{\bb_\phi, \
\phi=c\cdot{\phi^{\tau,\nu}_\pqr},\ c\in\k^\times\}$
 provides
a semiuniversal
 Poisson deformation  of
$\bb_{\phi^\tau}$, the coordinate ring of the corresponding
elliptic  singularity ~\eqref{pqr}.
\end{thm}

In the next section, we will state a `quantum analogue'
of Theorem \ref{delpezzoform} with Poisson algebras being replaced
by noncommutative algebras.

\begin{rem}
Observe  that
the family of {\em Poisson} algebras $\bb_\phi$, in part (iii),
depends on $\mu+1$ parameters. The reason is that, although
the underlying surface $\phi^\tau=0$ does not depend on the
extra-parameter $c\in\C^\times$, the corresponding
Poisson structure does.
\end{rem}

\begin{proof}[Proof of Theorem \ref{delpezzoform}]
Part (i) is a simple consequence of equations
$d=a+b+c$, and $p=d/a,q=d/b,r=d/c$, combined with formula
\eqref{mu}. Part (ii) is a well known classical result, cf. \cite{D}.

Next, let $S=\C^2\times S_p\times S_q\times
S_r\times\C^\times$. Here, the parameters $\tau,\nu$ form
coordinates in the first factor $\C^2$,
the affine linear spaces $S_p,S_q,S_r$ are spanned by the
corresponding polynomials in \eqref{pol}, and the parameter
$c$ gives a coordinate on the last factor $\C^\times$.
Thus, members of the  family  $\{\bb_\phi, \
\phi=c\cdot{\phi^{\tau,\nu}_\pqr},\ c\in\k^\times\}$
are parametrized by points of $S$.
Let $o\in S$ be the point corresponding to
vanishing parameters
$\nu,c,P,Q,R$, i.e., to the  Poisson algebra $\bb_{\phi^\tau}$.

To prove (iii), we must show that the classifying map
for our   family
of Poisson algebras induces a vector space isomorphism
$T_oS\iso \PH^2(\bb_{\phi^\tau})$. According to Theorem
\ref{compute}, cf. also the discussion at the beginning of this
subsection, we have $\PH^2(\bb_{\phi^\tau})=\J(\phi^\tau)\oplus\C$,
where the direct summand $\C$ corresponds to the 1-dimensional
space $\J^{(\pp)}(\phi^\tau).$ By part (i), we compute
$$\dim\PH^2(\bb_{\phi^\tau})=\dim\J(\phi^\tau)+1=
\mu+1=(p-1)+(q-1)+(r-1)+3=\dim S.
$$

It is easy to see that the map
$T_oS\to \J(\phi^\tau)\oplus\C$ we are interested in
is the natural map sending a polynomial $c\cdot{\phi^{\tau,\nu}_\pqr}$
to its resudue class in the Jacobi ring.
This map is injective. Hence, it must be an isomorphism,
due to the above equality of dimensions.
\end{proof}

\medskip{\large
\section{Main results}}\label{main_ssec}

\vskip 4mm
\subsection{Algebras defined by a potential.}\label{defs}
Let $V$ be a $\k$-vector space with basis
$x_1,\ldots,x_n$,
and let $ F=TV=\k\langle x_1,\ldots,x_n\rangle,
$ be the corresponding free tensor algebra.
The commutator quotient space $ F\nat= F/[F,F]$ is
 a $\k$-vector space with the natural basis formed by cyclic words in the
alphabet $x_1,\ldots,x_n.$
Elements of  $F\nat$ are referred to as {\em potentials}.

Let $\Phi\in F\nat$.
 For each $j=1,\ldots,n,$
one defines  $\partial_j \Phi\in F,$ the corresponding
partial derivative of the potential, by the formula
$$
\partial_j \Phi:=\sum_{\{s\;|\; i_s=j\}}x_{i_s+1}\,x_{i_s+2}\ldots  x_{i_r}\; x_{i_1}\,x_{i_2}\ldots
x_{i_s-1}\in \k\langle \xxx\rangle.
$$
We extend this definition to arbitrary elements
$\xi=(\xi_1,\ldots,\xi_n)\in
\C^n,$ by $\C$-linearity, i.e. we put $\pa_\xi\Phi:=\xi_1\cdot\partial_1
\Phi
+\ldots+\xi_n\cdot\partial_n
\Phi.\,$
This way, we get
a linear map $V^*\to TV,\,\xi\mto\pa_\xi\Phi.$

Many  interesting examples of Calabi-Yau algebras arise from
the following construction of algebras associated with a potential, cf. \cite{Gi}.
Given $\Phi\in F\nat$, introduce an  associative algebra
\beq{AF}
\A(\Phi):=F/\llb \pa\Phi\rrb=\C\langle  x_1,\ldots,x_n\rangle\big/\llb \pa_i
\Phi\rrb_{i=1,\ldots,n},
\eeq
 a quotient of  $F$ by the two-sided ideal
generated by all $n$ partial derivatives, $\pa_i
\Phi, \,i=1,\ldots,n,$ of the potential $\Phi$.

\subsection{Filtered setting}\label{fff}
Let each of the generators $x_k,\, k=1,\ldots,n,$
be assigned some {\em positive}  degree
$\deg x_k=d_k\geq 1.$ This makes $V$ a graded vector space,
with homogeneous basis $x_k,\,k=1,\ldots,n$. Thus, the tensor algebra
$F=TV=\C\langle  x_1,\ldots,x_n\rangle$ acquires a graded algebra structure
with respect to the induced
{\em total} grading $F=\bigoplus_{r\geq 0} F^{(r)} $
(not to be confused
with the standard grading on the tensor algebra; the latter
corresponds to the special case where $\deg x_k=1$ for all $k$).

One may  also view  $F$ as a filtered algebra, with an
increasing filtration $\k=F^{\leq 0}\sset F^{\leq 1}\sset\ldots,$
given by $F^{\leq r}=F^{(0)}\oplus\ldots\oplus F^{(r)}$.
The filtration, resp. grading, on $F$ gives
rise to a filtration $F\nat^{\leq k},\,k=0,1,\ldots,$ resp. grading
$F\nat=\oplus_r F\nat^{(r)}$, on the commutator quotient space $F\nat.$

The increasing filtration on $F$ induces a filtration
$\k=\A^{\leq 0}(\Phi)\sset \A^{\leq 1}(\Phi)\sset \A^{\leq
2}(\Phi)\sset\ldots,$  on
the quotient algebra $\A(\Phi)$.
In the special case where $\Phi$ is, in effect, homogeneous,   our algebra inherits a
 grading $\A(\Phi)=\bigoplus_{m\geq 0}\A^{(m)}(\Phi).$

Given a  filtered algebra $A$ with filtration by finite dimensional vector spaces, we write
$$\hilb(A):=
\sum\nolimits_{m\in\Z}\  \dim(\gr^{(m)}A)\cdot u^m \in \Z[[u]],
$$
for the  Hilbert-Poincar\'e series of the {\em associated graded} algebra
 $\gr A=\bigoplus_{m\geq 0}\gr^{(m)}A$.

An ascending filtration, resp. grading,
on $A$ induces a  filtration $\HH^\hdot_{\leq m}(A)$, resp. grading
$\HH^\hdot(A)=\bigoplus_{m\in\Z}\HH^\hdot_{(m)}(A),$ on each Hochschild
cohomology group.

 A  family of nonnegatively  filtered algebras is said to be a  semiuniversal
 filtered  family provided the associated graded algebras
   form a flat family
and,
moreover,  the  classifying map gives an isomorphism
$T_sS\iso \HH^2_{\leq0}(A_s)$ for all $s\in S.$
There is a similar definition in the case of graded algebras.

The above discussion also applies to filtered, resp. graded,
Poisson algebras and  Poisson cohomology.

\subsection{Quantization of the Poisson algebras $\aa_\phi$ and
$\bb_\phi$.}
\label{AABB}
In the  three sections below, we  are going to state
four theorems
which are  main results of the paper.
The proofs of these theorems will be given later,
mostly in \S\S\ref{deform_sec},\ref{coh_sec}.

 Fix a triple of integers  $0<a\leq b\leq c$ such that ${\mathsf{gcd}}(a,b,c)=1.$
We will be interested in (not necessarily commutative) algebras
with 3 generators. We put $F=\k\langle
x,y,z\rangle$ and view $F$ as a graded algebra
such that $\deg x=a,\,\deg y =b,\,\deg z=c.$
\smallskip

It will be convenient to introduce the following

\begin{defn}
An element $\Phi\in F\nat$ is called a {\em CY-potential}
provided   $\A(\Phi)$  is a Calabi-Yau algebra of {\bf dimension 3}.
\end{defn}

The basic example of a
homogeneous CY-potential of degree $d=a+b+c$ is  $\Phi=xyz-yxz
\in F\nat^{(d)}$.
In this case, one easily finds that
$\A(\Phi)=\mathbb C[x,y,z]$.

We will be mostly interested in general, not
necessarily homogeneous, potentials of degree $d= a+b+c$.

\begin{thm}\label{one}
 Let $(a,b,c)$ be a triple of positive integers
and $\Phi^{(d)}$ a  homogeneous CY-potential
of degree $d=a+b+c$. Then, for any  potential $\Phi'\in
 F\nat^{<d}$,  one has
\vskip 2pt

\vi \parbox[t]{153mm}{The sum $\,\Phi=\Phi^{(d)}+\Phi'$ is a
CY-potential, and for the corresponding filtered algebra, we have
$$
 \hilb\big(\A(\Phi)\big) \
=
\
1/(1-u^a)(1-u^b)(1-u^c),
$$
is the  Hilbert-Poincar\'e series of the graded algebra  $\k[x,y,z]$.}
\vskip 2pt

\vii  \parbox[t]{153mm}{There exists a nonscalar
central element $\Psi\in \A^{\leq d}(\Phi)$.}
\end{thm}

Theorem \ref{one} is proved in Subsection \ref{1and3}.

The equation in part (i) of the theorem shows that any
algebra of the form $\A(\Phi)$,
where $\Phi$ is a nonhomogeneous  potential such that
its leading term is a CY-potential of degree $a+b+c$, may be
thought of as a `noncommutative analogue' of the polynomial
algebra $\k[x,y,z]$. Further, a Calabi-Yau structure (i.e. a
noncommutative volume) on
the algebra may be thought of as a noncommutative deformation
of a unimodular Poisson structure on the polynomial
algebra. As we will see in \S\ref{3} below, any
such unimodular Poisson algebra must be of the
form $\aa_\phi$ for
an appropriate polynomial $\phi\in\C[x,y,z]$.
Moreover, the polynomial $\phi$ is necessarily a
{\em central} element for the Poisson structure.

This suggests to view a central element $\Psi\in \A(\Phi)$
 as a noncommutative analogue of the polynomial
$\phi$. Thus,  one may view any algebra of the form
\beq{B}
\B(\Phi,\Psi):=\A(\Phi)/\llb\Psi\rrb,
\qquad\Psi\in Z(\A(\Phi)),
\eeq
(a quotient of the CY algebra $\A(\Phi)$ by the two-sided ideal
generated by the central element $\Psi$),
as  a noncommutative analogue of  a Poisson algebra
of the form $\bb_\phi=\aa_\phi/(\phi)$.

\subsection{Noncommutative del Pezzo surfaces.} \label{NCDP}
 For the rest of  section \ref{main_ssec} we assume that  $(a,b,c)$ is one of the triples from
table \eqref{table} and recall the nonhomogeneous polynomials
$\dis\phi^{\tau,\nu}_\pqr$,
 of degree $d=a+b+c$ defined in \eqref{phiPQR}.
According to Theorem \ref{delpezzoform}(ii),
the algebra $\dis\bb_\phi,\ \phi=\phi^{\tau,\nu}_\pqr,$ gives
the coordinate ring of an affine del Pezzo surface.

One the other hand, Theorem \ref{one}(ii) insures  the
existence of nontrivial central elements in the noncommutative
algebra $\A(\Phi)$.
Therefore, it is natural to look for cyclic potentials $\Phi$
of the form similar to one given by formula \eqref{phiPQR}, and to view
the corresponding algebras $\B(\Phi,\Psi)$, in \eqref{B}, as  quantizations
of those  del Pezzo surfaces.

To implement this program, fix complex parameters  $t,c$.
To each triple $P\in \k[x],Q\in \k[y], R\in \k[z],$
of  polynomials given by formulas
~\eqref{pol},
of degrees $p,q,r,$ respectively,
 we associate
the following  potential
\beq{xyz} \Phi^{t,c}_\pqr=xyz-t\cd yxz+
c\cd\big[P(x)+Q(y)+R(z)\big]\ \in\ F\nat.
\eeq

Clearly, $\Phi^{t,c}_\pqr$ is a nonhomogeneous potential of degree
$d$.
The corresponding algebra
$\A(\Phi^{t,c}_\pqr)$ is a filtered algebra with
generators $x,y,z$, and the following 3 relations
\beq{XYZ}
xy-t\cd yx=c\cd\fra{dR(z)}{dz},\quad yz-t\cd
zy=c\cd\fra{dP(x)}{dx},\quad zx-t\cd xz=
c\cd\fra{dQ(y)}{dy}.
\eeq

We need the following
\begin{defn}\label{generic_def}
Let $X$ be an irreducible variety, thought of as a
variety of `parameters'.
We say that a property (P) holds for
generic parameters $x\in X$
if there exists a countable family, $\{Y_s\},$ of closed subvarieties
$Y_s\subsetneq X$, such that (P) holds for any
$x\in X\sminus (\cup_s Y_s).$
\end{defn}

Recall formula \eqref{dimJ} for the Milnor number $\mu=\dim\J(\phi)$
of an elliptic singularity.
The two theorems below are our main results about
noncommutative del Pezzo surfaces.

\begin{thm}\label{two} For generic parameters
$(t,c,P,Q,R)$, formula \eqref{xyz} gives a CY-potential,
and we have
\vskip 2pt

\vi The algebras
$\A(\Phi^{t,c}_\pqr)$, with  relations \eqref{XYZ},
form a {\em semiuniversal filtered}   family of
associative algebras  that depends on $\mu$ parameters.
\vskip 2pt

\vii The algebras of the form
$\B(\Phi^{t,c}_\pqr, \Psi),$ where $\Psi\in\A^{\leq d}(\Phi^{t,c}_\pqr)$ is a nonscalar central
element,
give a {\em semiuniversal} family of associative algebras that depends on $\mu+1$  parameters.
\end{thm}

A sketch of proof of Theorem \ref{two} is given in Subsection \ref{sketch2}.

Our presentation for the algebras $\B(\Phi,\Psi)$ in terms of generators
and relations is not completely explicit yet, since we have
not explicitly described central elements $\Psi$. This can be
done by a direct computation which has been carried out by
Eric Rains, see \S\ref{magma_sec} and \cite{R}.

Part (1) of the next Theorem gives
 a `parametrisation' of noncommutative
del Pezzo algebras similar to the one provided,
in the commutative (Poisson) case, by Theorem
\ref{delpezzoform}(iii).

\begin{thm}\label{three}
 For any {\em generic} homogeneous potential
$\Phi^{(d)}$ of degree $d=a+b+c$ and an arbitrary  potential $\Phi'\in
 F\nat^{<d}$,
the sum $\,\Phi=\Phi^{(d)}+\Phi'$ is a
CY-potential, and the following holds:
\vskip 1pt

\begin{enumerate}
\item
 There exists a  potential of the form $\Phi^{t,c}_\pqr,$ cf.
\eqref{xyz}, such that
 one has a filtered algebra isomorphism
$\dis\A(\Phi)\cong\A(\Phi^{t,c}_\pqr)$.
\vskip 2pt

\item The center of $\A(\Phi)$
is a free polynomial algebra $\k[\Psi]$
generated by
an element $\Psi\in\A^{\leq d}(\Phi),$ and one has  $\dis\gr
Z\big(\A(\Phi)\big)\cong
Z\big(\A(\Phi^{(d)})\big)$.
\end{enumerate}
\end{thm}

Theorem \ref{three} is proved in Subsection \ref{1and3}. 

\subsection{Noncommutative elliptic singularities} \label{NCES}
Let $(a,b,c)$ be one of the triples from table \eqref{table}.  In this subsection,
we are interested in the special case where the polynomials
$P,Q,R$, cf. \eqref{pol}, reduce to their leading terms. In such a case,
the corresponding potential $\Phi^{t,c}:=\Phi^{t,c}_\pqr$,
and the central element $\Psi\in \A^{\leq d}(\Phi^{t,c})$, both become
homogeneous elements of degree $\deg\Phi^{t,c}=\deg\Psi=a+b+c=d$.

Explicitly, we have, cf. also \S\ref{magma_sec} and \cite{R},

\beq{table2}
{\renewcommand{\arraystretch}{1.7}
\begin{tabular}{|c|c|c|}\hline
case&$\Phi^{t,c}\in F^{(d)}\nat$&$\Psi\in Z(\A^{(d)}(\Phi^{t,c}))$\\
\hline
$E_6$ &$
 xyz-t\cd yxz
+ c(\frac{x^3}{3}+\frac{y^3}{3}+\frac{z^3}{3})$&
$
c\cd y^3 + \frac{t^3 -   c^3}{c^3 +  1}(y
z   x+c\cd  z^3)
 - t
\cd  z   y   x$\\
\hline
$E_7$ & $ xyz-t\cd yxz
+c(\frac{x^4}{4}+\frac{y^4}{4}+\frac{z^2}{2})$&
$(t^2 + 1)
x  y  x  y  -\frac{t^4 +t^2 +1}{t^2 - c^4}(t\cd x
y^2  x+c^2\cd y^4) +t \cd y^2  x^2
$\\
\hline
$E_8$ & $xyz-t\cd yxz
+c(\frac{x^6}{6}+\frac{y^3}{3}+\frac{z^2}{2})$&
too long\\
\hline
\end{tabular}}
\eeq

Let $\chi(u)$ denote the rational function in the RHS of formula
\eqref{R}. Further, let $\ups\in \HH^3(\A(\Phi))$ denote the
image of $1\in\HH_0(\A(\Phi))$ under the isomorphism in \eqref{vdb}, resp.
$\De$ denote the BV-operator,
associated with a noncommutative volume
 on the CY algebra $\A(\Phi),$ cf. Definition
\ref{vdb_def}.

\begin{thm}\label{four} Let $(a,b,c)$ be as  in
table \eqref{table2}. Then, for any generic homogeneous potential $\Phi$
of
degree
$d=a+b+c$, one has
\vskip 2pt

\vi  There exists a  potential of the form $\Phi^{t,c},$  as in
table \eqref{table2},
 such that
 one has a graded algebra isomorphism
$\dis\A(\Phi)\cong\A(\Phi^{t,c})$.

\vii
Each group $\HH^k\big(\A(\Phi)\big),\, k\leq 3,$ is a free
$\k[\Psi]$-module with the
Hilbert-Poincar\'e series:
 $$
\hilb\big(\HH^k(\A(\Phi))\big)=
\begin{cases}
{\frac{1}{1-u^d}}_{_{}} &\oper{if}\en k=0,1;\\
\frac{1}{u^d}{\big[\frac{\chi(u)^{^{}}}{(_{_{}}1-u^d)_{_{}}}-1\big]}_{_{}}\quad &\oper{if}\en k=2;\\
{\frac{\chi(^{^{}}u)^{^{}}}{u^d(1-u^d)}}^{^{}} &\oper{if}\en k=3.
\end{cases}
$$

\viii The BV-operator kills $\ups$ and
induces the following bijections
$$
\De:\
\HH^3(\A(\Phi))/\k\cd\ups\iso \HH^2(\A(\Phi)),
\quad\oper{resp.}\quad
\De:\
\HH^1(\A(\Phi))\iso\HH^0(\A(\Phi)).
$$
\end{thm}

Theorem \ref{four} is proved in Subsection \ref{pf4}.

\begin{rems} (1)\en Part (ii) of the theorem is a
 generalization of a result of Van den Bergh \cite{VB2}.
The factor $u^d$ in denominators of the formulas
is due to the fact that $\deg\ups=-d$.
We recall also that any Calabi-Yau algebra of dimension 3 has no
Hochschild cohomology in degrees ~$>3$.

(2)\en For a result related to part (i) see also
\cite{BT}, Proposition 5.4.
\end{rems}

Associated with a nonzero homogeneous central element
$\Psi\in\A(\Phi)$, of degree $d$, there is the corresponding quotient algebra
$\B(\Phi,\Psi)$, cf. \eqref{B}, which inherits a graded algebra structure.
According to \cite{ATV} and \cite{St1},
the element $\Psi$ is not a zero divisor in $\A(\Phi);$ furthermore,
the algebra $\B(\Phi,\Psi)$ is a noetherian domain of  Gelfand-Kirillov dimension two.

Let $D^b({\B(\Phi,\Psi)})$ be the bounded derived category of
finitely generated graded left $\B(\Phi,\Psi)$-modules.
One also introduces $\Tails(\B(\Phi,\Psi))
\sset D^b({\B(\Phi,\Psi)}),$ a full triangulated
subcategory of tails, whose objects are
complexes with finite dimensional cohomology, cf. \cite{NVB}.

Recall next that, for any algebra of the form
$\B(\Phi,\Psi)$ as above,
there exists a triple
 $(\E,{\mathcal L},\sigma)$, where $\E$ is an elliptic curve,
${\mathcal L}$ is a positive line bundle on $\E$, and $\sigma$ is an automorphism of $\E$,
such that one has a graded algebra isomorphism, see  \cite{ATV}, \cite{St1},
\beq{proj}
\B(\Phi,\Psi)=\bigoplus\nolimits_{m\geq 0}\ \Gamma(\E, {\mathcal L}\o \sigma^*{\mathcal L}\o\ldots
\o(\sigma^{m-1})^*{\mathcal L}).
\eeq

The graded algebra on the right of \eqref{proj} is a $\sigma$-twisted
homogeneous coordinate ring of $\E$.
Therefore, the algebra $\B(\Phi,\Psi)$ may be thought of as a flat graded
noncommutative deformation of the
affine cone over the elliptic curve $\E$.

Let  $D^bCoh(\E)$ be the bounded
 derived category
of coherent sheaves
on $\E$.
According to a result due to Artin and Van den Bergh, \cite{AV},
there is a triangulated equivalence
\beq{artvdb}
D^bCoh(\E)\cong D^b({\B(\Phi,\Psi)})/\Tails(\B(\Phi,\Psi)).
\eeq

\subsection{Matrix factorizations on a noncommutative singularity}
\label{MATR}
Given a nonnegatively graded algebra $A$ and a central homogeneous element $\Psi\in
A$,
of degree $d>0$, one
 may introduce  $D_{\operatorname{gr}}(A,\Psi),$ a triangulated
category of graded matrix factorizations, see  \cite{Or1}.
An object of $D_{\operatorname{gr}}(A,\Psi)$ is a diagram
\beq{factor}
M=\left(\xymatrix{
M_+ \ar@/^/[rr]^<>(0.5){g}&& M_-\ar@/^/[ll]^<>(0.5){g'}
}\right)
\qquad g\ccirc g'= \Psi\cd\Id_{M_-},\en g'\ccirc g= \Psi\cd\Id_{M_+},
\eeq
where $M_+,M_-$ is a pair of finite rank
free graded $A$-modules and
$g,g'$ is a pair of graded  $A$-module morphisms
of degrees $0$ and $d$, respectively.

We take $A=\A(\Phi)$ and
apply a noncommutative version of results due to Orlov,
\cite{Or1},\cite{Or2}.
This way,
one obtains the following, cf. also \cite{KST}.

\begin{thm}\label{orlov_prop} \vi There is a triangulated equivalence
$$D^bCoh(\E)\cong D_{\operatorname{gr}}(\A(\Phi),\Psi).
$$

\vii Any maximal Cohen-Macaulay graded $\B(\Phi,\Psi)$-module
has a 2-periodic free  graded $\A(\Phi)$-module resolution.
\end{thm}

\begin{proof}[Sketch of proof of Theorem \ref{orlov_prop}]
It is known that $\A(\Phi)$, being a graded Calabi-Yau algebra,
is automatically a Gorenstein, Artin-Schelter regular algebra
of dimension 3,
see \cite{BT}, \cite{ATV}. Further, the central element
$\Psi$ is not a zero divisor in $\A(\Phi)$, by construction.
It follows that the quotient $\B(\Phi,\Psi)=\A(\Phi)/\llb\Psi\rrb$
is an Auslander-Gorenstein algebra of dimension 2, by \cite{Le}.

Let $\Perf(\B(\Phi,\Psi))$
denote the full triangulated subcategory in  $D^b({\B(\Phi,\Psi)})$
of {\em perfect complexes},
i.e. of bounded complexes of  free graded
 left $\B(\Phi,\Psi)$-modules of finite rank. Following Orlov, \cite{Or1},
one introduces a quotient category
$D^\text{sing}_\text{gr}(\B(\Phi,\Psi)):=D^b({\B(\Phi,\Psi)}/\Perf(\B(\Phi,\Psi))$,
 the  triangulated  category of a homogeneous singularity.

An immediate
generalization of \cite{Or1}, Theorem 3.9, yields
the following result

\begin{prop}\label{or1} Let $A=\oplus_{j\geq 0} A_j$ be a
graded noetherian  algebra
with $A_0=\k$. Assume that $A$ is
 Gorenstein, Artin-Schelter regular algebra of dimension $n$.
Let $\Psi\in A_n$ be a homogeneous central
element, which is not a zero divisor. Then,
there is a triangulated equivalence
$$
D^{\operatorname{sing}}_{\operatorname{gr}}(A/\llb\Psi\rrb)\cong
D_{\operatorname{gr}}(A,\Psi).
$$
\end{prop}

\begin{proof}
The proof of this result is based on the fact
that $A/\llb\Psi\rrb$ is
 an Auslander-Gorenstein algebra of dimension $n-1$,
by \cite{Le}. This insures that an analogue of \cite{Or1}, Proposition
1.23, holds in our present noncommutative setting. The rest of the proof
 of \cite{Or1}, Theorem 3.9 then goes through, and Proposition \ref{or1}
 follows.
\end{proof}
\smallskip

Next, we apply \cite{Or2}, Theorem 2.5, to the algebra $\B(\Phi,\Psi)$.
This way, we obtain  a triangulated equivalence
\beq{quot}
D^\text{sing}_\text{gr}(\B(\Phi,\Psi))\cong
D^b({\B(\Phi,\Psi)})/\Tails(\B(\Phi,\Psi)).
\eeq

On the other hand, applying the equivalence of Artin and Van den Bergh,
 \eqref{artvdb},
 and using the isomorphism in \eqref{proj},
we deduce that the quotient category on the right of
\eqref{quot} is equivalent to  $D^bCoh(\E)$.
This, combined with Prposition \ref{or1}, yields part (i) of
Theorem
\ref{orlov_prop}.

The proof of part (ii) is similar to the proof  of  the corresponding well-known
result in commutative algebra, due to D. Eisenbud \cite{Ei}.
\end{proof}

\begin{examp}\label{mat_examp} One of the simplest examples is the
case of a cubic curve $\E_\tau\sset\PP^2=\PP(\k^3),$ with homogeneous equation
of the form, cf. \eqref{pqr},
\beq{phip}
\psi^\tau(x,y,z):=x^3+y^3+z^3 +\tau\cd xyz,\qquad \tau\in\C^*.
\eeq

Motivated by \cite{ATV} and \cite{LPP}, for any point  $u\in\PP^2,$
with homogeneous coordinates $(\al,\beta,\gamma),$
one associates the following $3\times 3$-matrix $D,$ as well as the corresponding
adjoint $D^\natural$, the matrix formed by the
$2\times 2$-minors of $D$,
$$
D:=\left(
\begin{matrix}
\alpha x & \beta z & \gamma y\\
\gamma z &\alpha y & \beta x\\
\beta y&\gamma x&\alpha z
\end{matrix}\right),\quad
D^\natural=\left(
\begin{matrix}
\alpha^2 yz -\beta \gamma x^2&\gamma ^2xy-\alpha \beta z^2& \beta ^2xy-\alpha \gamma y^2\\
\beta ^2xy-\alpha \gamma z^2& \alpha ^2yz-\beta \gamma y^2& \gamma ^2yz-\alpha \beta x^2\\
\gamma ^2xz-\alpha \beta y^2&\beta ^2yz-\alpha \gamma x^2&\alpha ^2xy-\beta \gamma z^2
\end{matrix}\right).
$$

We have  an identity
$D\cd D^\natural=D^\natural\cd D=\det D\cd\Id$.
Assume  that  $\alpha, \beta, \gamma$ are all nonzero and put
$D':=-\frac{1}{\alpha\beta\gamma}\cdot D^\natural.$ Thus,
we obtain an equation
$D\cd D'=D'\cd D=-\frac{\det D}{\alpha\beta\gamma}\cd\Id$.

 Further, from the definition of $D$ one computes
$$\det D=(\alpha^3+\beta^3+\gamma^3)xyz-\alpha\beta\gamma(x^3+y^3+z^3).$$

Therefore, assuming that the triple $(\alpha, \beta, \gamma)$
is such that
 $\alpha^3+\beta^3+\gamma^3=\tau\cdot\alpha\beta\gamma$
we may write, $\det D=-\alpha\beta\gamma\cd\psi^\tau.$
We deduce that whenever $\psi^\tau(\alpha,\beta,\gamma)=0$ holds, cf.
\eqref{phip},
one has $D\cd D'=\psi^\tau\cd\Id=D'\cd D.$
This way, we have constructed
a
family of graded matrix factorizations
\beq{fact3}M_{u}=\left(\xymatrix{
\C[x,y,z]^{\oplus 3} \ar@/^/[rr]^<>(0.5){D}&& \C[x,y,z]^{\oplus 3}\ar@/^/[ll]^<>(0.5){D'}
}\right)\in  D_{\operatorname{gr}}(\C[x,y,z],\psi^\tau),
\qquad{u}\in\E_\tau,
\eeq
parametrized by the points ${u}=(\alpha,\beta,\gamma)\in \E_\tau$
with nonvanishing coordinates.
\hfill$\lozenge$\end{examp}

There is an important class of  {\em point modules}
over the algebra $\A(\Phi)$ introduced in \cite{ATV}. A point module
has a grading $P=\bigoplus_{k\geq 0} P^{(k)}$ such that
$P^{(0)}=\C$ and $\dim P^{(k)}\leq 1$ for any $k$. Given an integer
$r>0$, we let $P_{\leq r}:= P/\bigoplus_{k> r} P^{(k)}$ denote
the $r$-truncation of $P$.

 Following  \cite{ATV}, one proves that any
point module $P$ is annihilated by $\Psi$, hence, may
be viewed as a $\B(\Phi,\Psi)$-module. Further, it is not difficult to
show
that
there exists $r>0$ such that the map $P \mto P_{\leq r}$
assigning to a point module its $r$-truncation gives
a bijection between the moduli spaces of point modules
and $r$-truncated point modules, respectively. Let $r_o$ be the minimal
such $r$.

We expect that Example \ref{mat_examp} can be generalized to a noncommutative
setting.
Specifically, let $\Phi$ be a homogeneous CY potential of degree
$d=a+b+c$, and let $\Psi\in\A(\Phi)$ be a homogeneous central
element of degree $d$.

\begin{conj}\label{matr_conj} To any  point module $P$ over the algebra
$\B(\Phi,\Psi)$ one can associate naturally a matrix
factorization $M(P)=(M_+,M_-),$ as in \eqref{factor},
where $\rk M_\pm=\dim P_{\leq r_o}$.
\end{conj}

In the
$E_6$-case, one has
 $r_o=d-1=2$ and $\dim P_{\leq r_o}=d=3.$
Moreover,
 it was shown in  \cite{ATV} that
point modules are parametrized by the points of the corresponding
elliptic curve $\E$. In that case,
our conjectural matrix factorisation $M(P)$ should reduce
to \eqref{fact3}, where $u\in \E$ stands for the parameter of
the point module $P$.

\medskip{\large
\section{Three-dimensional Poisson structures}\label{3}}
\vskip 4mm

\subsection{}\label{not}
Given a (not necessarily
smooth)
finitely generated commutative $\k$-algebra  $ A $, write $\Om^1 A $ for the $ A $-module of
K\"ahler differentials of $ A $, and let
 $\Om^\hd A :=\La^\hd_ A (\Om^1 A )$
be the graded commutative algebra of differential
forms, equipped with the de Rham differential $\dr$.
 For each $p=1,2,\ldots,$
we also have
$\X^p A  =\Hom_ A (\Om^p A, A ),$ the space of skew $p$-polyderivations
$ A \wedge_\k\ldots\wedge_\k A \to A .$

Set  $\X^0 A  := A =: \Om^0 A .$
The graded
space $\X^\hd A  :=\bigoplus_{p\geq0}\X^p A  $ has a natural structure of
Gerstenhaber algebra with respect to the Schouten bracket
 $[-,-]:\X^p A  \times \X^q A
\to\X^{p+q-1} A  .$
Associated with a polyderivation $\eta\in \X^p A  $,
there is  a Lie derivative operator
$L_\eta:  \Om^\hdot A \to\Om^{\hdot-p+1} A ,$
resp. contraction  operator
$i_\eta:\Om^\hdot  A \to\Om^{\hdot-p} A .$
These operators make $\Om^\hdot A $ a Gerstenhaber
$\X^\hd A  $-module.

Let $ A =\k[\mm]$ be the coordinate ring of a
{\em smooth} affine variety $\mm,$ with tangent
bundle $T_\mm,$  resp. cotangent bundle $T^*_\mm$.
Then we have canonical isomorphisms $\X^\hd A  =
\Gamma(\mm, \wedge^\hd T_\mm),$ resp.
$\Om^\hdot  A =\Gamma(\mm, \wedge^\hd T^*_\mm).$
We will also use the notation $\X^\hd(\mm),$ resp.
$\Om^\hdot(\mm)$,
for these spaces.

\subsection{Unimodular Poisson structures.}
Any Poisson bracket $\{-,-\}: \ A\times A\to A$ on a (not necessarily smooth) finitely generated
commutative algebra $A$ determines (and is determined by)
a bivector $\pi\in\X^2 A, $ via the  formula
\beq{bra}
\{f,g\}:=\lan \dr f\wedge \dr g,\pi\ran,\qquad\forall f,g\in A.
\eeq
The Jacobi identity for the bracket $\{-,-\}$ is  equivalent
to the equation   $[\pi,\pi]=0,$ in $\X^3 A  $.

Associated with any $f\in  A, $ there is
a Hamiltonian derivation $\xi_f:=\{f,-\}\in \X^1 A  ;$
it is easy to check that $\xi_f=[\pi,f].$

Let $\mm$ be a smooth affine variety  of dimension $n$,
with a trivial canonical bundle.
Let
 $\vol\in \Om^n(\mm)$ be
 a nowhere vanishing volume $n$-form.
Contraction with $\vol$ yields an isomorphism
\beq{XOm}
\X^p(\mm)\ \iso\
\Om^{n-p}(\mm),\quad \eta\mto i_\eta\vol,\qquad
p=0,\ldots,n.
\eeq

A Poisson bracket  on the algebra $A=\k[\mm]$ is said to be {\em unimodular}
provided the divergence (with respect to the volume $\vol$)
of any Hamiltonian vector field vanishes, i.e,
for any $ f\in\k[\mm]$, we have $\div(\xi_f)=0.$
This means that
 the volume-form  is preserved by
 the Hamiltonian flow generated by the vector field $\xi_f$.

One has the following  standard result.

\begin{lem}\label{unimod}
Given  an arbitrary  bivector $\pi\in\X^2(\mm) $, on a {\em 3-dimensional} smooth
  variety $\mm$,
let $\al:=i_\pi\vol$, a 1-form. Then,  we have\smallskip

\vi \en\parbox[t]{140mm}{The condition that $\pi$ be
a Poisson bivector is equivalent to the equation
$\al\wedge \dr \al=0.$}
\smallskip

\vii$\en\dis
\pi \text{ gives a unimodular Poisson bracket}
\quad\Longleftrightarrow\quad
L_\pi\vol=0
\quad\Longleftrightarrow\quad
\dr \al=0.
$
\end{lem}

\begin{proof} For any $\eta\in \X^p(\mm) ,$
one has $i_{[\pi,\eta]}=[L_\pi, i_\eta],$
where $[-,-]$ stands for the {\em super}-commutator.
 Further, using Cartan's identity
$L_\pi:=i_\pi \dr  - \dr i_\pi$  we get
$$
i_{[\pi,\eta]}=
i_\pi \dr  i_\eta -\dr  i_\pi i_\eta - (-1)^p i_\eta i_\pi \dr
+(-1)^p i_\eta \dr  i_\pi.
$$

We take $p=2$ and apply the operations on each side of the
 identity  to the 3-form $\vol$. Clearly, one has
$\dr\ \vol=0$ and also $i_\pi i_\eta\vol =i_\eta i_\pi\vol=i_\eta \al=0$. Hence, we find
\beq{pieta}
i_{[\pi,\eta]}\vol=i_\pi \dr  i_\eta\vol +i_\eta \dr \al +\dr  i_\pi \al=
\langle\pi, \dr  i_\eta\vol\rangle+\langle\eta,\dr \al\rangle.
\eeq

Now let $\eta=\pi$ and
let $\ups$ be the 3-vector inverse to $\vol$.
Then, we have $\pi=i_\al\ups .$ So,
$\langle\pi, \dr \al\rangle=\langle i_\al\ups ,\dr \al\rangle=
\langle\ups, \al\wedge \dr \al\rangle.$ Hence, we obtain
$i_{[\pi,\pi]}\vol=2\langle\ups, \al\wedge \dr \al\rangle.$
Thus, we see that $[\pi,\pi]=0$ holds if and only if
we have $\al\wedge \dr \al=0.$
This yields part (i) since the pairing in \eqref{bra}
gives a Poisson bracket
if and only if   one has $[\pi,\pi]=0$.

There is also an alternate more geometric proof of (i)
as follows.
A bivector $\pi$  gives a Poisson structure on $\mm$ if and only if
$[\pi,\pi]=0$, which holds
if and only if   the  distribution in $T_\mm $ (of generic  rank 2) spanned by $\pi$
is integrable.  For $\al=i_\pi\vol$,
the same distribution may
be alternatively described as the  distribution defined by the kernels
of the 1-form $\al$. The latter distribution is integrable if and only if
$\al$ satisfies  Frobenius integrability condition:
$ \al\wedge \dr \al=0.$

The unimodularity property in part (ii) is equivalent to the equation
\beq{pffor}
0=\div(\xi_f)\cd\vol=L_{\xi_f}\vol=\dr (i_{\xi_f}\vol),\qquad\forall
f\in\k[\mm].
\eeq

We have
\beq{xi}
\xi_f=i_{\dr f}\pi=i_{\dr f}(i_\al\ups)=i_{\dr f\wedge\al}\ups .
\eeq
Therefore, we get  $i_{\xi_f}\vol=\dr f\wedge\al,$ hence
$\dr (i_{\xi_f}\vol)=-\dr f\wedge \dr \al.$ We see that \eqref{pffor}
amounts to the equation $\dr f\wedge \dr \al=0,$ for
any regular function $f$. This holds if and only if   we have $\dr \al=0.$
\end{proof}
\subsection{}  Fix a smooth 3-dimensional manifold
with a nowhere vanishing volume form $\vol\in\Om^3(\mm)$
and a regular function $\phi\in \k[\mm]$.

 Associated with
$\dr \phi$, an  {\em exact}  1-form,
one has a bivector $\pi\in\X^2(\mm) $ such that
$i_\pi\vol=\dr \phi.$
By Lemma
\ref{unimod}, this bivector gives rise
to a  unimodular Poisson  bracket  $\{-,-\}_\phi$, on $\k[\mm].$
Explicitly, the bracket  is determined by the equation
\beq{fg}\{f,g\}_\phi\cdot\vol=\dr \phi\wedge \dr f\wedge \dr g,
\qquad\forall f,g\in\k[\mm].
\eeq

We now specialize to the case where $\mm=\k^3,$ is a vector
space with coordinates $x,y,z,$ and
$\vol=\dr x\wedge dy\wedge dz$ is the standard  volume
form.

\begin{cor}\label{alg_closed} Let  $\{-,-\}$ be an unimodular polynomial Poisson
structure on $\k[x,y,z]$. Then,

\vi There exists a  polynomial $\phi\in\k[x,y,z]$,
such that the   Poisson bracket of linear functions
is given by formula \eqref{formula}.

\vii We have $\k[\phi]\sset\zz(\k[x,y,z]).$ If   the  Poisson bracket is
  nonzero then, any element $ f \in \zz(\k[x,y,z])$
  is algebraic over the subalgebra
$\k[\phi],$ i.e.  there exists
a nonzero polynomial $P\in\k[t_1,t_2]$ such that one has $P(\phi, f )=0$.
\end{cor}

\begin{proof} Recall that any polynomial closed 1-form on $\k^3$ is exact. Hence,
any  unimodular Poisson bracket on
the algebra $\k[\mm]=\k[x,y,z]$ is of the form
\eqref{fg}, for some
  polynomial  function $\phi\in \k[\mm].$
The corresponding  Poisson bivector
$\pi$ is given by
\beq{Vol}
\pi=i_{\dr\phi}\ups=
\mbox{$
\frac{\partial \phi}{\partial x}\cdot  \frac{\pa}{\pa y}\wedge
\frac{\pa}{\pa z}+\frac{\partial \phi}{\partial y}\cdot   \frac{\pa}{\pa z}\wedge \frac{\pa}{\pa x}
+\frac{\partial \phi}{\partial z} \cdot  \frac{\pa}{\pa x}\wedge
\frac{\pa}{\pa y}$},
\quad\text{where}\quad\ups:=\mbox{$\frac{\pa}{\pa x}\wedge \frac{\pa}{\pa y}\wedge
\frac{\pa}{\pa z}$}.
\eeq

Part (i), and the inclusion in part (ii) follow.

Next, let $ f  \in\k[x,y,z]$ be  such that
$\k(\phi, f )$, the field of rational functions generated
by the polynomials $\phi$ and $ f $, has
transcendence degree $= 2$ over $\k$. Then, there exists
a point $u\in\k^3$ such that $\dr\phi|_u$ and $\dr f |_u$ are linearly
independent covectors.

Now, formula \eqref{fg} shows that $ f $ is a central element with
respect to the Poisson bracket if and only if   one has
$\dr\phi\wedge \dr f =0.$ Hence, for $ f \in\zz(\k[x,y,z]),$
the covectors   $\dr\phi|_u$ and $\dr f |_u$ must be proportional,
and part (ii) follows.
\end{proof}
\begin{rem}\label{centr_rem}  For a polynomial $\phi$ such that the ring
$\C[\phi]$ is algebraically closed in $\C[x,y,z]$,
Corollary \ref{alg_closed}(ii)  yields
$\k[\phi]=\zz(\k[x,y,z]).$
 This condition holds for instance
for any irreducible  polynomial, cf.
\cite{P}, Proposition 4.2, for a similar result in a special case.
\end{rem}

\medskip{\large\section{Poisson (co)homology}}
\vskip 4mm\subsection{}\label{hp_subsec}
Poisson homology $\PH_\idot(A)$, resp. cohomology $\PH^\hdot(A)$, of a Poisson algebra $A$
is defined as the homology  of the total complex associated with
a  double complex,
$DP_{\idot,\idot}( A )=\La^\hd_ A  (D_\idot\Om^1 A ) ,$ resp.
  $DP^{\hd,\hd}( A )=\Hom_ A (DP_{\idot,\idot}( A ) , A ),$
cf. \cite[Appendix]{GK}.
Here, $D_\idot\Om^1 A  $ denotes the {\em cotangent complex}
of $ A $, the latter being viewed as a commutative associative algebra, cf.
\cite{GK}, formula (A.4).

The bi-graded space $DP^{\hd,\hd}( A )$ comes equipped
with a natural Gerstenhaber (i.e. graded Poisson) algebra structure, of bi-degree
$(0,-1)$ that gives rise to a  Gerstenhaber algebra structure on  $\PH^\hdot(A)$,
see \cite{GK}. Also,
$\PH^0(A)=\zz(A)$ and for each $j=0,1,\ldots,$  the group $\PH^j(A)$, resp.
$\PH_j(A)$,
comes equipped with a natural $\zz(A)$-module structure.

Let $\pi\in\X^2 A  $ be the bi-derivation associated with the
Poisson bracket, cf. \eqref{bra}.
The Lie derivative  $L_\pi: \Om^\hd A\to
\Om^{\hd-1} A,$
resp.   $L_\pi: \X^\hd A   \to\X^{\hd+1} A  $,  makes
the graded space $\Om^\hd A$, resp. $\X^\hd A,$
 a complex called homological, resp.
 cohomological,
 {\em Lichnerowicz-Koszul-Brylinski complex} (LKB-complex) of the
 Poisson algebra $ A $,
cf. \cite{Br}.

The canonical projection
$D_\idot\Om^1 A  \onto \Om^1 A  $ induces a map
$DP_{\idot,\idot}( A )\onto DP_{0,\idot}( A )=\Om^\hd A $,
and also the dual
map $\X^\hd A  \to DP^{\hd,\hd}( A ).$
These maps provide morphisms between the  LKB-
and Poisson cohomology complexes, respectively.  Furthermore,
unlike the case of the Hochschild complex,
the map  $\X^\hd A  \to DP^{\hd,\hd}( A )$
 turns out to be
a  DG Gerstenhaber algebra morphism.

If the scheme
$\Spec A $ is  smooth then the  projection $D_\idot\Om^1 A  \onto \Om^1 A  $
is a quasi-isomorphism.
It follows that each of the morphisms
$DP_{\idot,\idot}( A )\to\Om^\hd A $,
and  $\X^\hd A  \to DP^{\hd,\hd}( A ),$ is
 a quasi-isomorphism as well.
In that case, Poisson (co)homology of $ A $ may be computed
via the corresponding  LKB  complex, that is,
 one has,  cf. ~\cite{GK},
$$\PH_\idot( A )=H^\hd(\Om^\hd A ,\ L_\pi),\quad
\text{resp.}\quad
\PH^\hdot( A )=H^\hdot(\X^\hd A  ,\  L_\pi).
$$

Observe that the de Rham differential $\dr: \Om^\hd A\to \Om^{\hd+1}A$
anti-commutes with the operator $L_\pi$, hence, induces
a well defined operator $\dr: \PH_\idot(A)\to \PH_{\idot+1}(A),$
on Poisson homology,
cf. \cite{Xu}.

Assume next that  $\Spec A$ is a  manifold of pure dimension $n$,
equipped with a nowhere vanishing volume form $\vol\in \Om^n A$.
Define a differential  $\delta: \X^\hd A\to \X^{\hd-1} A,$
 by transporting  the de Rham differential
$\dr:\Om^\hd A\to\Om^{\hd+1}A$ via the isomorphism
$\X^\hd A\iso \Om^{n-\hd}A,$ cf. \eqref{XOm}.
Then, by \cite{Xu}, Proposition 4.5 and Theorem 4.8,
we have

\begin{prop}\label{weinstein} Let $\Spec A$ be  smooth of pure dimension $n$.
 For any  unimodular Poisson bivector $\pi\in\X^2 A$, one has

\vi The isomorphism in \eqref{XOm}
intertwines the $L_\pi$-actions on polyvector fields and
on differential forms; it
 induces  a degree
reversing $\zz(A)$-module isomorphism $\dis\PH_\idot(A)\iso \PH^{n-\hd}(A).$

\vii The differential  $\delta$
anti-commutes with  $L_\pi$; it descends
to a well-defined BV-type differential $\delta:
\PH^\hd(A)\to \PH^{\hd-1}(A)$.\qed
\end{prop}

\subsection{Poisson homology of a complete intersection.}
Let
 $I\sset A$ be a  Poisson  ideal in a
Poisson algebra $A$, so we have  $\{I,A\}\sset I$.
We set $ B=A/I.$ Thus, $ B$ is  a Poisson algebra and $\Spec  B$
is a closed Poisson subscheme in $\Spec  A$.

The following is a Poisson analogue of a similar result known
in the case of Hochschild cohomology, cf. eg. \cite{LR}.

\begin{prop}\label{AJ} Assume that the Poisson scheme
 $\Spec A$ is  {\em smooth}
and, moreover,  the Poisson subscheme $\Spec  B$ is a locally complete
intersection in   $\Spec A$. Then, one has $L_\pi(I^n\cdot \Om^mA)
\sset I^{n+1}\cdot \Om^{m-1}A,$  for any
$m,n\geq 0,$ and
there is a direct sum decomposition
$$\PH_k( B)=\bigoplus_{0\leq 2j\leq k}
H^{k-2j}\big(I^j\cd\Om^\hd A/I^{j+1}\cd\Om^\hd A,
\ L_\pi\big),\qquad \forall k\geq 0.
$$
\end{prop}
\proof The first statement is verified by a direct computation.
 Further, the assumption that $\Spec  B$ be a  locally complete
intersection insures that $I/I^2$ is a projective $ B$-module,
and the cotangent complex of  $\Spec  B$ may be represented
by a length two complex of amplitude $[-1,0]$,
$$D_\idot\Om^1 B\
 \stackrel{\text{qis}}\simeq\
\big[I/I^2 \stackrel{d}\too  B\o \Om^1A \big].
$$

Hence,
Poisson double complex is quasi-isomorphic to a  double complex with
the following terms
\begin{align*}
\qquad DP_{p,q}(B)\
\stackrel{\text{qis}}\simeq \
\La^q_B (D_p\Om^1 B)
&=\La^q_B\big([I/I^2][1]\ \bigoplus\ B\o_A \Om^1A \big)&\\
&=\bigoplus_{0\leq j\leq q}\big([\Sym^j(I/I^2)][j]\o_A\Om^{q-j}A\big)&\\
&=\bigoplus_{0\leq j\leq q}\big(I^j\cd\Om^{q-j} A/I^{j+1}\cd\Om^{q-j}A\big)[j].&\qquad\qquad\qed
\end{align*}

\subsection{Poisson cohomology of a hypersurface.}\label{spectr}
Below, we will be mostly interested in Poisson {\em co}homology
of an algebra of the form $B:=\C[\mm]/(\phi)$
where $\mm$ is a smooth Poisson variety and $\phi\in  \C[\mm]$ a regular
function contained in the Poisson center. In that case,
one can give a slightly different
 description of  Poisson
(co)homology of the algebra $B$, which is  more explicit
than the one provided by Proposition \ref{AJ}.

Observe first that contraction with the 1-form $\dr\phi$
provides a  differential $i_{\dr\phi}: \X^\hd(\mm)\to\X^{\hd-1}(\mm),$
in the
corresponding Koszul
 complex.

\begin{rem}
The     Jacobi ring of $\phi$  may be
identified with $H^0(\X^\hdot (\mm),\, i_{\dr\phi}).$ The latter group
 is the only nontrivial cohomology group
of the  Koszul complex, provided $\phi$ has   isolated
singularities. This way, using
the Euler-Poincar\'e principle, one proves formula
\eqref{R}.
\end{rem}

Let  $\pi\in\X^2(\mm) $ be a Poisson
bivector.
\begin{lem}  For any $\phi\in\zz(\k[\mm])$,
the map $L_\pi$ is $\C[\phi]$-linear and it
anticommutes with $i_{\dr\phi}$.
\end{lem}
\begin{proof} In general, for any function $f$ and
a bivector $\pi$, one has the following standard identity
\beq{LL}
L_\pi\ccirc i_{\dr f}+ i_{\dr f}\ccirc L_\pi
=L_{i_{\dr f}\pi}=\xi_f.
\eeq

Now, the function $\phi$ is central with respect
to the Poisson bracket given by a bivector $\pi$
if and only if one has $i_{\dr\phi}\pi=0$.
In that case the maps  $L_\pi$ and $i_{\dr\phi}$
anticommute. The $\C[\phi]$-linearity
statement is clear.
\end{proof}

\begin{prop}\label{spectralka}
Let  $\phi\in\zz(\k[\mm])$ be a central
regular function on $\mm$. Assume that $\phi$ has   only
isolated
singularities and that there exists a vector field
$\eu\in\X^1(\mm)$ such that one has $\eu(\phi)=c\cdot\phi$,
where $c$ is a nonzero constant.

 Then, for the Poisson cohomology of the algebra ${B}_\phi:=\k[\mm]/(\phi)$,
there is a  convergent first quadrant
spectral sequence $\dis E^{p,q}_2\ \Rightarrow\ \gr_p \PH^{p+q}({B}_\phi),$
with $E_1$-term of the form

{\footnotesize
\beq{spsequence}
\xymatrix{
\dots&&\dots&&\dots&&\dots&&\dots\\
\X^4({B}_\phi)\ar[u]^<>(0.5){L_\pi}&&
0\ar[u]^<>(0.5){L_\pi}&&0\ar[u]^<>(0.5){L_\pi}&&\J(\phi)\eu
\ar[u]^<>(0.5){L_\pi}&&\J(\phi)\ar[u]^<>(0.5){L_\pi}
\\
\X^3({B}_\phi)\ar[u]^<>(0.5){L_\pi}&&
0\ar[u]^<>(0.5){L_\pi}&&\J(\phi)\eu
\ar[u]^<>(0.5){L_\pi}&&\J(\phi)\ar[u]^<>(0.5){L_\pi}\\
\X^2({B}_\phi)\ar[u]^<>(0.5){L_\pi}&&
\J(\phi)\eu\ar[u]^<>(0.5){L_\pi}&&
\J(\phi)\ar[u]^<>(0.5){L_\pi}&&\\
\X^1({B}_\phi)\ar[u]^<>(0.5){L_\pi}&&
\J(\phi)\ar[u]^<>(0.5){L_\pi}&&&&\\
\X^0({B}_\phi)\ar[u]^<>(0.5){L_\pi}&&&&&&
}
\eeq}

\noindent
\!\!where the leftmost column
is the  LKB  complex of the Poisson algebra
${B}_\phi$,
and $\J(\phi)$ denotes the Jacobi ring of $\phi$, cf.
\S\ref{coh}.
\end{prop}

\begin{proof} Put $A:=\k[\mm]$ and
let $\ta^\hdot= A\o\k[\tau]/(\tau^2)$ denote a graded {\em
super}-commutative
algebra such that $ A$ is an even subalgebra placed in degree zero,
and  $\tau$ is an odd generator  of degree $-1$. We introduce
a differential $\nabla:\ta^\hdot\to\ta^{\hd+1}$,
which is defined as an  odd super-derivation, $\nabla=\phi\frac{\pa}{\pa \tau}$,
that annihilates the subalgebra $ A$ and
is such that $\nabla(\tau)=\phi$.
Clearly, one can view the resulting DG algebra as a 2-term
complex $A\stackrel{\phi\cdot}\too A$ with the differential
given by multiplication by the function $\phi$.
Therefore, we have $H_0(\ta)={B}_\phi$ and $H_j(\ta)=0$
for any $j\neq0$.

Next, we make $\ta$ a Poisson DG algebra by extending
the Poisson bracket $\{-,-\}$ on $A$
by $\tau$-linearity. This way, $\ta$ becomes a  Poisson  DG algebra
which is quasi-isomorphic to ${B}_\phi$. Thus, for the Poisson
cohomology, we have
$\PH^\hdot({B}_\phi)\cong \PH^\hdot(\ta),$
where the cohomology on the right-hand side denotes
the hyper-cohomology involving the differential
$\nabla$, on $\ta$.

It will be convenient to use  geometric language
and write $\ta=\k[Y],$ where $Y=\mm\times\k$ is a
 {\em smooth} affine
Poisson super-manifold  of super-dimension $(\dim\mm|1)$.
The corresponding Poisson cohomology  may be computed,
according to
general principles,
as a hyper-cohomology of the    LKB
double complex for the Poisson DG super-manifold $Y$.
This way, we deduce
$$\PH^\hdot( A)\cong \PH^\hdot(\ta)\cong
H^\hdot(\X^\hd\ta,\,\nabla+  L_\pi),
$$
where the differential  $\nabla$ is induced
by the same named differential on the DG algebra $\ta$ itself,
and  the differential
$ L_\pi$ comes from the Poisson
structure on $\k[\mm]$.

Let $t$ denote an {\em even} coordinate on
the total space, $T^*Y$,  dual to the
{\em odd} coordinate
$\tau$ on $Y$. Thus, we get $\X^\hd(Y)=\X^\hd(\mm) \o\k[t,\tau]/(\tau^2),$
where the variable $t$ is assigned grade degree $+2$.
With this notation,  the    LKB  complex
takes the form
\beq{form1}\left(\X^\hdot(\mm) \o\k[t,\tau]/(\tau^2),\
\phi\cd \mbox{$\frac{\pa}{\pa
\tau}$}+t\cd i_{\dr\phi}+  L_\pi\right).
\eeq

Let $T_\phi:=T_\mm |_{\phi\inv(0)}$ denote
the restriction of the tangent bundle of $\mm$
to the (not necessarily smooth)
hyper-surface $\phi\inv(0)\sset\mm$.
Thus,
$T_\phi$ is a vector bundle on $\phi\inv(0)$
of rank $\dim\mm$, and we let
 $\La^\hdot_\phi:=\Gamma(\phi\inv(0),
\La^\hdot T_\phi)\cong{B}_\phi\o_{ A}\X^\hd A$ be the corresponding
exterior algebra viewed as a graded algebra such that
the space ${B}_\phi\o_{ A}\X^1 A$ is placed in  degree $+1$.

Restriction to  $\phi\inv(0)$ combined with the
specialization $\tau\mto0,$ gives a natural algebra projection
\beq{pr}
\mathsf{pr}: \
\Big(\X^\hd(\mm) \o\k[t,\tau]/(\tau^2),\
\phi\cd \mbox{$\frac{\pa}{\pa\tau}$}+t\cd i_{\dr\phi}+  L_\pi\Big)
\onto \Big(\La^\hdot_\phi\o\k^\hdot[t],\ t\cd i_{\dr\phi}\Big).
\eeq
It is easy to see that the differential
in \eqref{form1} descends to the  differential $t\cdot i_{\dr\phi},$
on $\La^\hdot_\phi\o\k[t]$. Moreover,
 the map $\mathsf{pr}$
is a quasi-isomorphism of DG algebras.

Thus, we conclude that the Poisson cohomology of ${B}_\phi$
may be computed as hyper-cohomology of the
DG algebra represented by the following  {\em mixed complex}:
{\small\beq{mixed}
\xymatrix{
\ldots&\ldots&\ldots\\
\La^2_\phi\ar[u]^<>(0.5){L_\pi}\ar[r]^<>(0.5){i_{\dr\phi}}&
\La^1_\phi\ar[u]^<>(0.5){L_\pi}\ar[r]^<>(0.5){i_{\dr\phi}}&
\La^0_\phi\ar[u]^<>(0.5){L_\pi}\\
\La^1_\phi\ar[u]^<>(0.5){L_\pi}\ar[r]^<>(0.5){i_{\dr\phi}}&
\La^0_\phi\ar[u]^<>(0.5){L_\pi}&\\
\La^0_\phi\ar[u]^<>(0.5){L_\pi}&&
}
\eeq}

We view
\eqref{mixed} as a bicomplex, $K$, with two differentials,
$i_{\dr\phi}$ and $L_\pi$.
Associated with this bicomplex, there is a
standard first quadrant spectral sequence $(E^{p,q}_r, d_r)$ such that
$E_1=H^\hdot(K, i_{\dr\phi})$ and the differential $d_1$ is induced by
$L_\pi$.

We first analyze the horizontal differential $i_{\dr\phi}$.
Let $\Lambda^{(q)}:\ \La^q_\phi\to\ldots\La^1_\phi\to\La^0_\phi$
denote the complex in the $q$-th row of diagram \eqref{mixed}
where, for any $j=0,1,\ldots,q,$
the term $\La^j_\phi$ is placed in degree $j$.
The $E_1$ page of the spectral sequence
of the bicomplex \eqref{mixed} takes the required form
\eqref{spsequence} thanks to the sublemma below.
\end{proof}
\begin{sublem}\label{xker} \vi We have $H^q(\Lambda^{(q)},\,i_{\dr\phi})=
\X^q{B}_\phi.$

\vii The complex $\Lambda^{(q)}$
is acyclic in all degrees $j\neq 0,1,$ or $q$.
Moreover,  $H^0(\Lambda^{(q)},\,i_{\dr\phi})=\J(\phi)$ and, if  $q>1$,
also we have
$H^1(\Lambda^{(q)},\,i_{\dr\phi})=\J(\phi).$
\end{sublem}
\begin{proof}[Proof of Sublemma]
To prove (i), recall that, in general,
 for any complete intersection
$\mathcal N\sset \mm$ where $\mm$ is smooth, each poly-derivation of
the algebra $\k[\mathcal N]$ is induced by a section
of the vector bundle $\La^p T_\mm|_{\mathcal N}.$
We take $ \mathcal N=\mm_\phi:=\phi\inv(0).$
It follows that each poly-derivation  $\theta: {B}_\phi\wedge\ldots\wedge{B}_\phi\to
{B}_\phi$ comes from a section
$s\in \La^p_\phi=\Gamma(\mm_\phi,\La^p T_\mm|_{\mm_\phi})$.
An extension of $s$ to a section $\wt s\in \Gamma(\mm,\La^p T_\mm)$
gives a
 poly-derivation $\wt\theta:  A\wedge\ldots\wedge A\to
 A,$ such that one has
$\wt\theta(\phi,a_1,\ldots,a_{p-1})\in\phi\cdot A,$
for  any $a_1,\ldots,a_{p-1}\in  A$.
In geometric language, the latter
condition  translates into the equation
$i_{\dr\phi}s=0$, for the original section $s$.
This proves (i).

Assume now that the function $\phi$ has an
isolated singularity. Then, the complex
$$
\X^{(q)}:\
\X^q(\mm)\stackrel{i_{\dr\phi}}\too\X^{q-1}(\mm)
\stackrel{i_{\dr\phi}}\too\ldots\too
\X^1(\mm)\stackrel{i_{\dr\phi}}\too
\X^0(\mm)
$$
is exact everywhere except possibly the leftmost and  rightmost
terms.  Furthermore,
the cokernel at the
rightmost term equals $\J(\phi).$

By definition, we have a short exact sequence of complexes
$0\to \X^{(q)}\to\X^{(q)}\to\Lambda^{(q)}\to0,$
where the morphism $\X^{(q)}\to\X^{(q)}$ is given by
multiplication by the function $\phi$.
 From the corresponding long exact sequence of cohomology,
we deduce that $H^j(\Lambda^{(q)},i_{\dr\phi})=0$
unless
$j\neq 0,1,q$. Moreover, since $H^0(\X^{(q)},i_{\dr\phi})=
\J(\phi)$,
 the final part of the long exact sequence reads
\beq{chase}0=H^1(\X^{(q)},i_{\dr\phi})\to
H^1(\Lambda^{(q)},i_{\dr\phi})\to \J(\phi)
\stackrel{\phi\cdot}\to\J(\phi)\to
H^0(\Lambda^{(q)},i_{\dr\phi})\to0.
\eeq

Now, by our assumptions, we have
$\phi\in\C\cdot i_{\dr\phi}\eu$.  Therefore,
the image of $\phi$ in the Jacobi ring
$\J(\phi)$ vanishes. Thus, the map
$\J(\phi)\to\J(\phi)$ induced by multiplication
by $\phi$ is equal to zero. This, combined with
the exact sequence \eqref{chase}, yields
part (ii) of the sublemma.
In addition, an easy diagram chase shows
that the preimage of the
element $1\in\J(\phi)$
under the isomorphism $H^1(\Lambda^{(q)},i_{\dr\phi})\iso\J(\phi)$,
cf. \eqref{chase}, corresponds to the class of the vector
field $\eu$.
\end{proof}

\subsection{Poisson cohomology of the algebra $\aa_\phi$.}
\label{pcoh}
We now specialize to the setting
of \S\ref{coh}. Thus, let $0<a\leq b\leq c$ be a triple of integers
with $\mathsf{gcd}(a,b,c)=1$.
Write $\eu$ for the corresponding Euler vector field
\eqref{eu}, on $\mm=\k^3,$ and  $\ups$ for the standard
3-vector, see \eqref{Vol}.

Given a polynomial homogeneous
  $\phi\in \k[\mm]=\k[x,y,z]$, write
$\pi:=i_{\dr\phi}\ups,$ cf.  \eqref{Vol},
and
let $\aa_\phi$ denote the
corresponding Poisson algebra, cf. Definition \ref{aabb}.

A. Pichereau has  found all   Poisson cohomology groups
of the algebra $\aa_\phi$ explicitly,
see \cite{P}. To state some of her results set $\mu:=\dim\J(\phi)$ and
choose homogeneous elements $1, f_1,\ldots, f_{\mu-1}\in \k[x,y,z]$
such that their residue classes modulo the Jacobi ideal form
a $\k$-basis of the vector space $\J(\phi).$
View the elements $i_{\dr f_k}\ups\in\X^2(\C^3),\, k=1,\ldots,\mu-1,$
 as elements of the LKB-complex
for the algebra $\aa_\phi$.

\begin{prop}[Pichereau]\label{HPA} For  any homogeneous
 polynomial  $\phi$ with an isolated singularity of degree
 $a+b+c$, Poisson cohomology of $\aa_\phi$ vanishes in degrees $\geq
 4$
and, one has
\vskip 2pt

\vi \parbox[t]{152mm}{We have $\PH^0(\aa_\phi)=\k[\phi].$ Furthermore,
the group
 $\PH^1(\aa_\phi)=\k[\phi]\eu$ is a rank 1 free
 $\k[\phi]$-module with generator ~$\eu$.}

\vii \parbox[t]{152mm}{The group $ \PH^3(\aa_\phi)$
is a rank $\mu$ free
 $\k[\phi]$-module with basis $\ups,\,
 f_1\cd\ups,\ldots,f_{\mu-1}\cd\ups$,
resp. $\PH^2(\aa_\phi)$ is a  free
 $\k[\phi]$-module with basis
 $i_{\dr f_1}\ups,\ldots,i_{\dr f_{\mu-1}}\ups,\,\pi.$\qed}
\end{prop}

\subsection{Poisson cohomology
of a  quasi-homogeneous singularity.}\label{Poihom} Pichereau has also computed cohomology groups
 of the
 LKB  complex for the (singular) Poisson algebra $\bb_\phi$
associated with
 a quasi-homogeneous polynomial  $\phi\in \k[x,y,z]$ of an {\em arbitrary} weight $d>0$.
Specifically,
she shows that $\X^p\bb_\phi=0,$ for all $p>2$. Furthermore,
the cohomology  of the LKB-differential
$L_\pi$ is as follows, see \cite{P},
\beq{piphi}
H^0(\X^\hd \bb_\phi)=\k,\quad
H^1(\X^\hd \bb_\phi)=\J^{(\pp)}(\phi)\eu,\quad
H^2(\X^\hd \bb_\phi)\cong\J^{(\pp)}(\phi)\pi,
\eeq
where
$\J^{(\pp)}(\phi)$ is viewed, in the notation of the previous subsection,  as the span
of the basis elements $f_j$ with $\deg f_j=\pp$.

\begin{proof}[Proof of Theorem  \ref{compute}]
We begin with  part (ii) of the theorem.
Our  Poisson bivector has the form
$\pi=i_{\dr\phi}\Upsilon$.
Therefore, for any $f\in \aa_\phi$ we have
$L_\pi(f)=\xi_f=i_{\dr f\wedge \dr \phi}\Upsilon
\sset i_{\dr \phi}(\X^2(\aa_\phi)).$
It follows that, for any $p\geq 1,$ the vertical differential
$L_\pi: E_1^{p,p}=\J(\phi)\to E_1^{p,p+1}=\J(\phi)\eu,
$ in  the spectral sequence
\eqref{spsequence}, vanishes. Thus, the $E_2$ page of the spectral sequence
reads

\beq{dspe}
\begin{matrix}
\dots&&\dots&&\dots&&\dots&&\dots\\
H^4(\X^\hdot\bb_\phi, L_\pi)&&
0&&0&&\J(\phi)\eu
&&\J(\phi)
\\
H^3(\X^\hdot\bb_\phi, L_\pi)&&
0&&\J(\phi)\eu
&&\J(\phi)\\
H^2(\X^\hdot\bb_\phi, L_\pi)&&
\J(\phi)\eu&&
\J(\phi)&&\\
H^1(\X^\hdot\bb_\phi, L_\pi)&&
\J(\phi)&&&&\\
H^0(\X^\hdot\bb_\phi, L_\pi)&&&&&&
\end{matrix}
\eeq

Here, the cohomology in the leftmost column is
provided by formula \eqref{piphi}, hence vanish in degrees
$>2$. Thus, we see
that all  differentials
$d_r: E_r^{p,q}\to E_r^{p-r+1,q+r},\,r\geq3,$
have zero range,
and the statement of part (ii) follows.

The case of Hochschild cohomology is quite similar.
Write $\HH^\hdot(-)=\HH^\hdot(-,b),$ where we have
explicitly indicated  the Hochschild differential $b$. Then,
using the notation $\nabla=\phi\frac{\pa}{\pa \tau}$,
we get
$$\HH^\hdot(\bb_\phi,\,b)=\HH^\hdot(\ta,\,\nabla+b)=
\HH^\hdot(\k[Y],\,\nabla+b)=H^\hdot(\X^\hd(Y),\,\nabla),
$$
where the last isomorphism is due to
the Hochschild-Kostant-Rosenberg
theorem applied to the smooth super-scheme $Y$.

One can now repeat the argument in the proof of Proposition \ref{spectralka}
and replace the complex $(\X^\hdot(Y),\nabla),$
by a smaller complex $(\La^\hdot_\phi\o\k[t], t\cdot i_{\dr\phi}),$
which is quasi-isomorphic to it, cf. \eqref{pr}.
This way, we see that the  Hochschild cohomology of the algebra
$\bb_\phi$ may be computed as hyper-cohomology of the
complex similar to
\eqref{mixed}, where the vertical differential $L_\pi$ is replaced
by zero. This yields part (i).
\end{proof}

\begin{rem} \vi Theorem \ref{compute}
shows that  Poisson cohomology groups of the algebra $\bb_\phi$
are  nonzero in all  degrees $\geq 2$, in particular,
these  groups are not the same as the cohomology groups  of the  LKB  complex,
cf. \eqref{piphi}. That agrees with the fact that
 the surface  $\mm_\phi=\Spec \bb_\phi$ has a singularity.

\vii Let $f\in \aa_\phi$.
For any $p=1,2,\ldots,$
the image of the element $f$ in $\J(\phi)$ gives
a class in $E_2^{p,p}$, cf. \eqref{dspe}.
An explicit lift of that class to
a $2p$-cocycle in the total complex
associated with the corresponding
bicomplex \eqref{mixed} is provided
by the element $f+i_{\dr f}\Upsilon\in
\La^0_\phi\oplus\La^2_\phi.$ Indeed,
we have  $L_\pi f=\xi_f=i_{\dr f\wedge\dr\phi}\Upsilon$,
cf. \eqref{xi}. Further, using \eqref{LL}, the fact that
 $L_\pi \Upsilon=0$ and unimodularity of the
Poisson bivector $\pi$ (Lemma \ref{unimod}), we get
$L_\pi i_{\dr f}\Upsilon=L_{i_{\dr f}\pi}\Upsilon=
L_{\xi_f}\Upsilon=0$.
Thus, we compute
$$(L_\pi+i_{\dr \phi})(f+i_{\dr f}\Upsilon)=
L_\pi f+L_\pi i_{\dr f}\Upsilon+i_{\dr \phi}i_{\dr f}\Upsilon
=i_{\dr f\wedge\dr\phi}\Upsilon+0+
i_{\dr\phi\wedge\dr f}\Upsilon=0.
$$

Similarly, for a homogeneous function $f\in\aa_\phi,$
of degree $\deg f=k$, the element
$(\varpi-k) f\Upsilon+f\eu\in\La^3_\phi\oplus \La^1_\phi$ gives,
for each $p=1,2,\ldots$, a
$(2p+1)$-cocycle in   the total complex
associated with the corresponding
bicomplex \eqref{mixed}. To see this, one may use the
identity $\pi\wedge\eu=
\deg\phi\cdot\Upsilon$,
to obtain the following equation, see
\cite[formula (27)]{P},
$$L_\pi(f\eu)=(k-\varpi)f\pi-\deg\phi\cd\phi\cd
i_{\dr f}\Upsilon=(k-\varpi)f\pi\
\operatorname{\mathrm{mod}}\ (\phi).$$

Further, we have $i_{\dr \phi}\eu=\eu(\phi)=(a+b+c)\phi$.
Thus, we find that, in $\La^\hdot_\phi=\X^\hdot(\mm)/(\phi)$,
one has
\begin{multline*}
(L_\pi+i_{\dr \phi})\big[(\varpi-k)f\Upsilon+f\eu\big]=(\varpi-k)
f i_{\dr \phi}\Upsilon+
L_\pi(f\eu)+f\eu(\phi)\\=
(\varpi-k)
f\pi+(k-\varpi)f\pi\
\operatorname{\mathrm{mod}}\ (\phi)
=0.
\end{multline*}
\end{rem}

\medskip{\large
\section{Classification results}\label{class_sec}}
\vskip 4mm

\subsection{Proof of Proposition \ref{alggeo}}\label{class}
Assume that the curve is not rational.
Let $a\le b\le c$.

If all the degrees are equal, then they
are equal $1$, and $\deg \phi\le 3$. In this case, the statement is
classical (the $E_6$ case).

Now assume that the degrees are not equal to each other.
In this case the leading power of $z$ is $\le 2$.
If this power were 1, the curve would be rational, so
it is $2$. Consider two cases.

Case 1. $a<b=c$. In this case $z^2$ comes together with
$zy$ and $y^2$, so for generic coefficients, by making a linear change
of $y,z$, we can kill $z^2$ and $y^2$, so the leading
term in $z$ will be linear. This shows that the curve is
rational, contradiction.

Case 2. $a\le b<c$. Then the leading term in
$z$ is $z^2$. so we get $2c\le a+b+c$, hence $c\le
a+b$. After a change of variable the equation
of the curve can be written as $z^2=g(x,y)$, where $g$ is
homogeneous of degree $2c\le 2(a+b)$. Consider two cases.

Case 2a. $a=b$. In this case, $g$ has degree $3$ or $4$.
If the degree of $g$ is $3$, then $a=b=2,c=3$, and this
the curve is rational, because the point $(x,y,z)$ is equivalent
to $(x,y,-z)$ in the weighted projective space. Thus it remains
to consider the case when $\deg(g)=4$, and thus
$a=b=1,c=2$. In this case, for generic coefficients we do
get an elliptic curve (the $E_7$ case).

Case 2b. $a<b$. Then the terms that can be present in
$g(x,y)$ are $y^3$ and terms that contain $y$ in power $\le 2$.
Thus for the curve not to be rational, the term $y^3$ must be present.
So $2c=3b$, and thus $b\le 2a$. If $b=2a$, then
$a=1,b=2,c=3,\d=6$, and for generic coefficients
we indeed get an elliptic curve (the $E_8$ case). On the other hand,
if $b<2a$, then $g$ cannot contain quadratic terms in $y$
(the only possible quadratic terms are $y^2,xy^2,x^2y^2$, and none
of them have the right degree). The only linear term in $y$
that can occur in $g$ is $x^3y$, in which case the curve is given
by $z^2=y^3+x^3y$, in weighted projective space of weights
$(4,6,9)$. This curve is rational, because the point
$(x,y,z)$ is equivalent to $(x,y,-z)$. Otherwise,
the curve is $z^2=y^3+x^p$, $4\le p\le 5$, in weighted projective space with
weights $(6,2p,3p)$. If $p=5$, the curve is rational since
$(x,y,z)$ is equivalent to $(x,y,-z)$. If $p=4$, the curve is
given by the equation $z^2=y^3+x^4$, with weights $(3,4,6)$,
and the curve is rational since $(x,y,z)$ is equivalent to
$(x,\varepsilon y,z)$, where $\varepsilon$ is a cubic root of
unity.
\qed

\subsection{Proof of Theorem \ref{three}(1)}\label{proof_one}
Let $Y'$ be the space of all non-homogeneous potentials of degree $a+b+c$,
and $Y$ be the space of all nonhomogeneous commutative polynomials of that
degree. Let $G'$ be the group of degree preserving automorphisms
of $\mathbb C\langle x,y,z\rangle$. Then we have the following exact sequence of
$G'$-modules:
$$
0\to U\to Y'\to Y\to 0,
$$
where $U$ is a 1-dimensional representation spanned by
$xyz-yxz$ in the $E_6$ case, and a 2-dimensional representation
spanned by $xyz-yxz$ and $xyxy-x^2y^2$ in the $E_7$ and $E_8$ cases.

Also, let $G$ be the group of degree preserving automorphisms
of $\mathbb C[x,y,z]$. We have an exact sequence
$$
1\to H\to G'\to G\to 1,
$$
where $H=1$ in the $E_6$ case, and $H=\mathbb G_a$
consisting of elements $x\to x,y\to y, z\to z+b(xy-yx)$
in the $E_7$ and $E_8$ cases. It is easy to see that a generic element
of $U$ is equivalent under $H$ to $\gamma(xyz-yxz)$.
Thus to prove the theorem, it suffices to show that
the expressions $xyz+c\cdot[P+Q+R]$, cf. \eqref{phiPQR}, give normal forms of generic
elements in $Y$ under the action of $G$. But this is a classical
fact from the theory of del Pezzo surfaces, see \cite{D}.
\qed

\medskip{\large
\section{Calabi-Yau deformations}\label{deform_sec}}
\vskip 4mm

\subsection{The dg algebra $\D(\Phi)$.}
Let $F=\k\langle
x_1,\ldots,x_n\rangle$ be  a free  algebra on $n$
 homogeneous generators
$x_1,\ldots,x_n$, where $\deg x_i>0$ for all $i=\nnn$. We  view
$F$ either as a graded or as
a filtered algebra, as in \S\ref{fff}.
We shall refer to the  grading on $F$ as
a {\em weight grading}.

Associated with any potential $\Phi\in F\nat$,
we have introduced in \cite{Gi}, \S1.4, a free graded  associative algebra
$\D(\Phi)=\oplus_{r\geq 0}\D(\Phi)_r $,  with $2n+1$
 homogeneous generators $x_1,\ldots,x_n, y_1,\ldots,y_n, t.$
We have
$\D(\Phi)_0=F$. The algebra $\D(\Phi)$ comes equipped
with a  differential $\pa: \D(\Phi)_\idot\to\D(\Phi)_{\idot-1}$
such that one has
$H^0(\D(\Phi))=F/\pa(\D_1(\Phi))=\A(\Phi).$

In the case where $\Phi$ is a homogeneous potential of degree $d> \text{max}\{\deg x_i, \,i=1,\ldots,n\},$
 there is an  additional weight grading on $\D(\Phi)$ such that the
generators $y_1,\ldots,y_n, t,$  are assigned degrees $\deg y_i:=d-\deg x_i,$ and
 $\deg t:=d$. This way, multiplication by elements of $\D(\Phi)_0$
makes each component
 $\D(\Phi)_r$  a {\em graded} left $F$-module
 $\D(\Phi)_r=\oplus_{s>0}\D^{(s)}(\Phi)_r $,
where the $s$-grading denotes the weight grading.

The precise definition of the dg algebra $\D(\Phi)$
is not essential for us at the
moment.  The important points are  the following 4 properties
\begin{align}
\bullet\en &{\text{\parbox[t]{140mm}{The differential on $\D(\Phi)$ and
the weight grading are determined by the
potential $\Phi$, while the  algebra structure and the $r$-grading do not
involve the potential;}}}\label{p1}\\
\bullet\en &{\text{\parbox[t]{140mm}{For each $r=0,1,\ldots,$
the homogeneous component $\D(\Phi)_r$  is a
free  $F$-module, moreover,
if $\Phi$ is homogeneous, then
 we have
$\dim_\k\big(\D^{(s)}(\Phi)_r\big)<\infty,\,\forall s\geq 0$;}}}\label{p2}\\
\bullet\en &{\text{\parbox[t]{140mm}{If $H^j(\D(\Phi))=0$ for all $j>0,$ then $\Phi$ is a
CY-potential;}}}\label{p3}
\\
\bullet\en &{\text{\parbox[t]{140mm}{If  $\Phi$ is a {\em homogeneous}
CY-potential then the differential $\pa$ preserves the weight grading on
$\D(\Phi)$; moreover, we have
 $H^j(\D(\Phi))=0$ for all $j>0,$
i.e.  the converse to \eqref{p3}
holds as well.}}}\label{p4}
\end{align}

Here, \eqref{p1}-\eqref{p2} are immediate from the definition
of $\D(\Phi)$, while \eqref{p3}-\eqref{p4} follow
from \cite{Gi}, Theorem 5.3.1, which is one of the
main results of that paper.
\medskip

For each $i=1,\ldots,n$, write $d_i:=\deg x_i >0,$
and let  $\Aut F$ denote the group of degree
preserving automorphisms of the algebra $F$.
Given $d\geq 3,$ let ${\rm CY_3}(d,d_1,\ldots,d_n)\sset F\nat^{(d)},$ be the set
 of homogeneous CY-potentials of some fixed degree
$d\geq 3$.

\begin{lem}\label{dense}
\vi The set ${\rm CY_3}(d,d_1,\ldots,d_n)$ is $\Aut F$-stable,
moreover, it
is an intersection of at most countable family of
Zariski open (possibly empty) subsets in $F\nat^{(d)}.$

\vii For all $\Phi\in {\rm CY_3}(d,d_1,\ldots,d_n)$,
the algebras $\A(\Phi)$ have the same
Hilbert-Poincar\'e series.
\end{lem}

\begin{proof} For any  $\Phi\in  F\nat^{(d)},$
we  may split the differential $\pa$ on
the dga  $\D(\Phi)$ into
components $\pa^\Phi_{r,s}  :\D^{(s)}(\Phi)_r\to  \D^{(s)}(\Phi)_{r-1},$
where each $\D^{(s)}(\Phi)_r$ is a finite dimensional vector space,
by \eqref{p2}. Since $\pa^2=0,$ for any $r,s\geq 0,$ one has
$\dim\Image\pa^\Phi_{r+1,s}\leq\dim\Ker\pa^\Phi_{r,s}$.

According to property \eqref{p4}, $\Phi$ is a CY potential
iff the dga  $\D(\Phi)$ has no nonzero cohomology in positive degrees.
Thus, we have
\beq{rank}
{\rm
CY_3}(d,d_1,\ldots,d_n)=\left\{
\Phi\in  F\nat^{(d)}\en\Big|\en
\begin{array}{c}
{\small\dim\Image\pa^\Phi_{r+1,s}\geq\dim\Ker\pa^\Phi_{r,s},\en\text{and}}\\
{\small\pa^\Phi_{1,s}\en\text{has maximal rank,}\en \forall r>0,s\geq0}.
\end{array}\right\}.
\eeq

The set on the right is clearly
 an intersection of a countable family of
Zariski open  subsets in $F\nat^{(d)}.$ Part (i) follows.
Part (ii) is \cite[Proposition 5.4.7]{Gi}.
\end{proof}

\subsection{Deformation setup}\label{setup}
In this and the following subsection, we develop a formalism 
that will be used in the
proofs of our main results.

Given a vector space $V$, we write
$V[[\h]]$ for the space of formal power series in an indeterminate $\h$ with
coefficients in $V$. In particular, we have
 $\k[[\hbar]]$,   the  ring of formal power series.
A  $\k[[\hbar]]$-module is said to
be {\em topologically free} if it is isomorphic to
a module of the form $V[[\h]]$, where $V$ is a $\k$-vector space.
Such a module is clearly a flat  $\k[[\hbar]]$-module,  complete
in  the $\h$-adic topology.
vvv
Let $\dis K=\oplus_{r\geq 0}K_r$
be a complex of  topologically free $\k[[\hbar]]$-modules, equipped with
a
$\k[[\hbar]]$-linear differential $\dis d: K_\idot\to
K_{\idot-1}.$
 Put $\dis \overline{K}:=K/\hbar\cdot
K$. This is a complex of $\k$-vector spaces,
with induced differential $\overline{d}:
 \overline{K}_\idot\to  \overline{K}_{\idot-1}.$

We recall  the following
standard result.

\begin{lem}\label{claim_complex} If the complex
$(\overline{K}, \overline{d})$ is acyclic in positive degrees then, we have

\vi The complex
$(K, d)$ is acyclic in positive degrees;

\vii The cohomology group $H^0(K, d)$ is
a  flat $\k[[\hbar]]$-module;

\viii
 The projection $K\onto \overline{K}$ induces an isomorphism
\dok 30mm
H^0(K, d)/\hbar\cd
H^0(K, d)\;\iso\;
H^0(\overline{K}, \overline{d}).
\edok
\end{lem}

We will also use a graded analogue of the above lemma,
where the variable $\hbar$ is assigned grade degee $1$.
Thus, let $K$ be a complex
of graded $\k[\hbar]$-modules
$K_r=\oplus_{s\geq 0}\ K_r^s$,
with homogeneous, $\k[\hbar]$-linear differential
$d: K_r\to K_{r-1}$.
Put $\overline{K}:=K/\hbar K$,
resp. $K':=K/(\hbar-1)K$, and
let $\overline{d}$, resp. $d'$, be the induced differential
on $\overline{K}$, resp. on $K'$.
For each $r$, the grading on $K_r$ induces
a filtration on $K'_r$. 
Replacing each term $K_r$ by its completion
$\wh K_r:=\prod_{s\geq 0} K^s_r$ and applying Lemma
\ref{claim_complex} to the resulting complex
yields the  following elementary result.
\begin{lem}\label{grclaim} Assume, in the above setting,
that each $K_r$ is a free
graded  $\k[\hbar]$-module such that $\dim_\C  K_r^s<\infty$
for all $r,s$, and that $H^r(\overline{K}, \overline{d})=0$
for any $r>0$. Then, we have $H^r(K',d')=0$ for
all $r>0$. Furthermore, the natural map
 $\overline{K}\to\gr K'$
induces an isomorphism
$H^0(\overline{K}, \overline{d})
\iso\gr H^0(K',d').$\qed
\end{lem}

Below, it will be
necessary to work with $\k[[\hbar]]$-algebras, that is,
with associative algebras $B$ equipped with a {\em central}
algebra imbedding $\k[[\hbar]]\into B$.
A  $\k[[\hbar]]$-algebra $B$ which is
complete in the  $\hbar$-adic topology
will be referred to as an $\h$-{\em algebra}.
Abusing terminology, we call such an algebra  {\em flat} if it is
topologically free as a left (equivalently, right) $\k[[\hbar]]$-module.

We reserve the notation $F_\h$ for the $\h$-algebra $F\hh.$
We have a canonical isomorphism of
free $\kh$-modules
$F_\hbar/[F_\hbar, F_\hbar]\cong F\nat[[\h]].$
This way, for any potential
\beq{expand2}\Phi=\Phi_0+\hbar\cdot\Phi_1+\hbar^2\cdot\Phi_2+\ldots
\in (F\hh)\nat=F\nat[[\h]],
\eeq
where $\Phi_j\in F\nat,$ one may define the following
$\h$-algebras

$$
\A_\hbar(\Phi):=F_\h\big/
\llb \pa_i\Phi\rrb_{i=1,\ldots,n},
\quad\text{resp.}
\quad \D_\hbar(\Phi)=\oplus_{r\geq 0}\D_\h(\Phi)_r.
$$
Here, $\D_\hbar(\Phi)$ is a
dg $\h$-algebra with $\kh$-linear  differential,
of degree $-1$, moreover, $\D_\h(\Phi)_0=F_\h,$ and we have
$\A_\h(\Phi)=H^0(\D_\h(\Phi)).$
There are natural `$\hbar$-analogues' of properties \eqref{p1}-\eqref{p4}.

\begin{cor}\label{dcor} Let $\Phi$ be a potential as in \eqref{expand2}.
Then,
we have

\vi For each  $r=0,1,\ldots,$ the component $\D_\h(\Phi)_r$ is a  free
 $F_\h$-module.

In the case where
all $\Phi_j$ are homogeneous of the same degree $d$,
the  homogeneous component
$\D^{(s)}(\Phi)_r $  is a finite rank free $\kh$-module,
for any $s\geq 0$.

\vii Reduction modulo $\hbar$ induces a dg algebra isomorphism
$\dis\D(\Phi_0)\iso\D_\h(\Phi)/\h\cd\D_\h(\Phi)\,
$
which, in the homogeneous case,
 is compatible with the weight gradings on each side.
\end{cor}
\begin{proof} Part (i) follows from an `$\h$-analogue'
of property \eqref{p2}; part (ii) follows from definitions.
\end{proof}

\subsection{Formal deformations of potentials}\label{def_sub}
For any vector space, resp. algebra, $C$ let $C\hb $ be
the vector space, resp. algebra,  of formal Laurent series with coefficients
in $C$. In particular, we put $\kk=\k\hb $,
the field of  Laurent series.

It is clear that $\qq F_\h=F\hb $. We also have
$F\nat[[\h]]\sset F\nat\hb =
[F\hb ]\nat.$ Therefore,
 any potential $\Phi\in F\nat[[\h]]$
may also be  viewed as a potential for the $\kk$-algebra $F\hb $.
Thus, one may view $\kk$ as a ground field and form a  $\kk$-algebra
$\A(\Phi)=F\hb /\llb\pa_j\Phi\rrb_{j=1,\ldots,n}$.
To emphasize the fact that the latter is a $\kk$-algebra,
we will write $\A_\kk(\Phi):=\A(\Phi)$.
There is an obvious $\kk$-algebra isomorphism
$\A_\kk(\Phi)=\qq\A_\h(\Phi).$

We begin with the following result which  says that being a CY-potential is
 an `open condition'.

\begin{prop} \label{open} Fix  a homogeneous  CY-potential
 $\Phi_0\in F^{(d)}\nat.$
\vskip 3pt

\vi  For any  (not necessarily homogeneous) element
$\Phi'\in F\nat\hh$, the sum  $\Phi=\Phi_0+\h\cdot\Phi'$  is a
 Calabi-Yau potential for the algebra $F\hb $.

Furthermore,
 $\A_\h(\Phi)$ is a {\em flat}
$\h$-algebra and the natural projection
yields an algebra isomorphism
$$\A(\Phi_0)\iso\A_\h(\Phi)/\h\cd\A_\h(\Phi).
$$

\vii For any  element
$\Phi'\in F\nat^{< d},$ the sum $\Phi=\Phi_0+\Phi'$
is a CY-potential for the algebra $F$.

Furthermore, the natural projection
yields a graded algebra isomorphism
$$\A(\Phi_0)\iso\gr\A(\Phi).
$$
\end{prop}
We remark that part (ii) of  Proposition \ref{open} is due to
 Berger and Taillefer, \cite{BT}; cf. also
\cite{Gi}, Corollary 5.4.4,
for an alternate approach.

\begin{proof}[Proof of Proposition \ref{open}]
To prove (i),  let $K:=\D_\h(\Phi)$.
Corollary \ref{dcor}(i) insures that the
assumptions of  Lemma \ref{claim_complex}
hold for $K$. It follows  from  property \eqref{p4}
and Lemma \ref{claim_complex}(i) 
that the dg algebra $\D_\h(\Phi)$ is acyclic in
positive degrees. Hence, the
dg algebra $\D_\kk(\Phi)=\qq\D_\h(\Phi)$ is acyclic in
positive degrees as well. Thus,  property \eqref{p3}
implies  that $\Phi$ is a Calabi-Yau potential.
Now, part (ii) of Lemma \ref{claim_complex}
insures that $\A_\h(\Phi)$  is a  flat
$\h$-algebra and part (iii) of Lemma \ref{claim_complex}
completes the proof of Proposition \ref{open}(i).

Now, we prove part (ii) of   Proposition \ref{open}
by an  argument involving Rees algebras that
 will be also  used  later in this section again.
Let $F^\hdot_\h:= F[\h]=\k[\h]\o F$.
We  assign the variable $\h$ degree $+1.$
This, together with  the grading $F=\oplus_s F^{(s)},$
makes $F^\hdot_\h$ a  graded $\k[\h]$-algebra,
the {\em Rees algebra} of $F$,
the latter being viewed as a {\em filtered} algebra.
 Thus, we have that $(F^\hdot_\h)\nat$
is a graded $\k[\h]$-module.

Next,  write a decomposition
$\Phi'=\Phi^{(d-1)}+\Phi^{(d-2)}+\ldots+\Phi^{(0)},$
 into homogeneous components $\Phi^{(r)}\in F\nat^{(r)},\,
 r=1,\ldots,d.$
Introduce a new {\em homogeneous}
potential (of degree $d$) for the graded algebra $F^\hdot_\h=F[\h]$ as follows
\beq{phiflat}
\Phi^\hbar:=\Phi_0+\hbar\cdot\Phi^{(d-1)}+\hbar^2\cdot\Phi^{(d-2)}+\ldots
+\h^d\cdot\Phi^{(0)}\in (F^\hdot_\h)\nat=F\nat[\h].
\eeq

One has a dg algebra 
$\D_\h^\hdot(\Phi^\hbar)=\oplus_{r,s\geq 0}\D_\h^{(s)}(\Phi^\hbar)_r,$
with differential
 $\D^{(\hdot)}_\h(\Phi^\hbar)_r
\to\D^{(\hdot)}_\h(\Phi^\hbar)_{r-1}$
defined in terms of the homogeneous potential $\Phi^\hbar$.
For each $r\geq 0,$ the
component $\D^{(\hdot)}_\h(\Phi^\hbar)_r$
is a free graded $\k[\hbar]$-module and the
differential is a morphism of graded $\k[\hbar]$-modules.
Further, we have 
dg algebra isomorphisms, cf. Corollary \ref{dcor}(ii):
\beq{ddaa}
\D_\h^\hdot(\Phi^\hbar)/(\h-1)\cdot\D_\h^\hdot(\Phi^\hbar)\iso\D(\Phi),
\quad\text{resp.}\quad
\D_\h^\hdot(\Phi^\hbar)/\h\cdot\D_\h^\hdot(\Phi^\hbar)\iso\D(\Phi_0).
\eeq

Here, the dg algebra on the right
 is acyclic
in positive cohomological degrees by \eqref{p4},
since $\Phi_0$ is a homogeneous CY potential.
Hence, the dg algebra on the left 
 is acyclic
in positive cohomological degrees, by Lemma \ref{grclaim}.
Also, from \eqref{ddaa}, we deduce 
$$
H^0\big(\D_\h^\hdot(\Phi^\hbar)/(\h-1)\cdot\D_\h^\hdot(\Phi^\hbar)\big)
\cong\A(\Phi),
\quad\text{resp.}\quad
H^0\big(\D_\h^\hdot(\Phi^\hbar)/\h\cdot\D_\h^\hdot(\Phi^\hbar)\big)
\cong\A(\Phi_0).
$$
Thus, the last statement of Lemma \ref{grclaim}
yields the algebra isomorphism
$\A(\Phi_0)\cong\gr\A(\Phi).$
\end{proof}

\subsection{The case: $n=3$}\label{3var}
We  put $F=\k\langle x,y,z\rangle$
and assign the generators $x,y,z$ positive weights $(a,b,c)$.
 Let $d:=a+b+c$.

First of all, we know that $\Phi_0:=xyz-yxz$ is a
CY-potential of degree $d$. In other words,
we have $\Phi_0\in {\rm CY_3}(d,a,b,c)$. We recall Definition
\ref{generic_def}, and deduce

\begin{cor}\label{open1}
\vi A generic homogeneous potential $\Phi\in F^{(d)}$
 is a CY-potential;
the Hilbert-Poincar\'e series of
the corresponding graded algebra $\A(\Phi)$ is equal to
that of the algebra $\C[x,y,z]$.

\vii For any $\Phi'\in F\nat^{<d}$,
the sum $\Phi=xyz-yxz+\Phi'$
is a CY-potential;
moreover, the natural projection
yields a graded algebra isomorphism
$$\mathbb C[x,y,z]\iso\gr\A(\Phi).
$$
\end{cor}
\begin{proof}
Part (ii) follows from  Proposition
\ref{open}(ii). Further, we observe that the set  ${\rm CY_3}(d,a,b,c)$ contains $\Phi_0$,
hence  is nonempty.
Therefore, part
(i) follows from  Lemma \ref{dense}.
\end{proof}

Recall that $\kk=\k(\!(\h)\!)$.
Since ${\rm CY_3}(d,a,b,c)\neq\emptyset$ for $d=a+b+c$, from Proposition
\ref{open}(i) we deduce

\begin{lem}\label{hphi}  For any element $\Phi'\in F\nat\hh,$ the sum
$\Phi=xyz-yxz+\h\cdot\Phi'$
is a CY-potential for the $\kk$-algebra $F\hb $.

Furthermore, the $\h$-algebra
$\A_\h(\Phi),$ with
relations
\begin{equation}\label{psi}
xy-yx=\hbar\cd \fra{\partial\Phi'}{\pa z},\quad yz-zy=\hbar\cd \fra{\partial\Phi'}{\pa x},\quad
zx-xz=\hbar\cd \fra{\partial\Phi'}{\pa y},
\end{equation}
is a flat  formal deformation of the polynomial algebra $\mathbb
C[x,y,z]$. \qed
\end{lem}

Reducing the flat deformation of the lemma
modulo $\hbar^2$, one obtains in a standard way
a Poisson bracket on $\mathbb C[x,y,z]$.
To describe this
Poisson bracket,  consider  the natural abelianization map
$$\C\langle x,y,z\rangle\nat\onto \C[x,y,z],\quad
f\mto f^{\oper{ab}}.
$$

Further, expand the potential in Lemma \ref{hphi} as a power series
in $\hbar$ and  write
\beq{expand}
\Phi= xyz-yxz +\hbar\cdot \Phi_1 +\hbar^2\cdot \Phi_2+\ldots,
\qquad
\Phi_j\in \C\langle x,y,z\rangle\nat.
\eeq

It is easy to show that the Poisson bracket on $\mathbb C[x,y,z]$
arising from the  flat deformation of Lemma \ref{hphi}
is given by formula \eqref{formula}; specifically, we have
\beq{defbra}
\{-,-\}=\{-,-\}_\phi
\quad\text{where}\quad\phi:=(\Phi_1)^{\oper{ab}}\in \C[x,y,z],
\eeq
the image  under
the abelianization map of
the degree 1 term in the $\hbar$-power series expansion of ~$\Phi$.

\medskip{\large
\section{From Poisson to Hochschild cohomology}\label{coh_sec}}
\vskip 4mm

\subsection{} We fix a triple of positive weights $(a,b,c)$.
 Put $F=\C\langle x,y,z\rangle$ and
assign the generators $x,y,z$ some positive weights $a,b,c$,
 respectively. This gives  the  ascending filtration
$F^{\leq m},\, m=0,1,\ldots,$
 on $F$, as in \S\ref{fff}. Further, we introduce a variable $\h$ of
{\em degree zero}
 and use the notation  $F^{\leq m}_\h:=(F^{\leq m})[[\h]]$, resp.
 $F^{\leq m}\nat[[\h]]$,
for the corresponding induced filtrations
on the $\h$-algebra $F_\h$, resp. on  $(F_\h)\nat=(F\nat)[[\h]].$
Thus, given a  potential $\Phi\in (F_\h)^{\leq m}\nat$,
we get a filtered $\h$-algebra $\A_\h(\Phi)$.
Note that these filtrations on $F_\h,\ (F_\h)\nat,$ and  $\A_\h(\Phi),$
are {\em not exhaustive},
rather, one has that
 e.g. $\dis\cup_{m\geq 0}\ \A^{\leq m}_\h(\Phi)$ is {\em dense} in
 $\A_\h(\Phi)$
with respect to  the $\h$-adic topology.

Now, put $d=a+b+c,$ and recall the notation
$\kk=\C((\hbar))$, resp. $\A_\kk(\Phi)=\kk\bigotimes_{\k[[\h]]}\A_\h(\Phi)$
and  Definition \ref{aabb}.

\begin{prop}\label{cen} \vi For any  potential
 $\Phi\in (F_\h)^{\leq d}\nat$ of the form \eqref{expand},
with $\Phi_0=xyz-yxz$, the  $\h$-algebra $\A_\h(\Phi)$
contains a central element
 $\Psi\in \A^{\leq d}_\h(\Phi)$
such that one has
$
\dis\Psi\ (\oper{mod}\h)=(\Phi_1)^{\oper{ab}}.
$
\vskip 4pt

\quad Assume, in addition, that  $(\Phi_1)^{\oper{ab}}$
is  a
homogeneous polynomial  of
degree $a+b+c$ with
 an isolated singularity.
Then we have:
\vskip 2pt

\vii There is a {\em bi-graded} $\kk$-algebra isomorphism
$$
\HH^\hdot(\A_\kk(\Phi))\cong
\kk\o\PH^\hdot(\aa_\phi)\quad\oper{where}\quad\phi:=(\Phi_1)^{\oper{ab}}\in\k[x,y,z].
$$

\viii The center of $\A_\h(\Phi)$ is
 $\dis Z(\A_\h(\Phi))=\k[\Psi][[\h]]$,
 a free topological $\h$-algebra in one generator, and $\dis
\HH^1(\A_\kk(\Phi))=\kk[\Psi]\Eu$, is a rank one free
$ \kk[\Psi]$-module generated by the Euler derivation.
\end{prop}

\begin{proof}[Proof of Proposition \ref{cen}$\mathsf{(i)}$]
Let $R=\k[u]$ be a graded polynomial algebra where the variable $u$ is
assigned grade degree $1$.
Below, we  consider $R$ as a ground ring,
and write $R[x,y,z]=\k[x,y,z,u]$, a polynomial
$R$-algebra.
Given a commutative $R$-algebra $A$ we use the notation
$\Om^\hd_RA,$ resp. $\X^\hd_RA,$ for the algebra of
{\em relative} differential forms
with respect to the subalgebra $R\sset A$, resp. $R$-linear
polyderivations of $A$.

Given a filtered algebra $B$ we write $\r B^\hd=\sum_{m\geq 0} B^{\leq
m}\cdot u^m$
for the corresponding  Rees algebra, a flat graded $R$-algebra.
Thus, associated with the filtered algebra
$F$, resp. $F_\h$, one has a  graded $R$-algebra
$\r F$, resp.  a  graded $R[[\h]]$-algebra
$\r F_\h$.

Now, fix a potential $\Phi=\sum \h^j\Phi_j\in (F_\h)^{\leq d}\nat,$ as in \eqref{expand},
and
for each $j=1,2,\ldots,$ write
$\Phi_j=\Phi_j^{(d)}+\Phi_j^{(d-1)}+\ldots+\Phi_j^{(0)},$
where $\Phi_j^{(m)}\in F\nat^{(m)}.$ We introduce a new
{\em homogeneous} potential of degree $d$ similar to the one in \eqref{phiflat}
(but
where  the role of $\h$  is now played by the variable $u$),
$$
\Phi^u=xyz-yxz+\sum_{j=1}^\infty\h^j\cd\Phi^u_j\ \in\
\r F\nat[[\h]],\quad\text{where}\quad
\Phi^u_j:=
\sum_{m=0}^d\ u^m\cd\Phi^{(d-m)}_j\in
\r^{(d)}F\nat.
$$

 Associated with the potential $\Phi$, resp. $\Phi^u,$ we have a filtered
$\h$-algebra $\A_\h(\Phi),$ resp.
 graded $R[[\h]]$-algebra
 $\A_\h(\Phi^u)$.
Clearly,
there is a natural graded algebra isomorphism
 $\r\A_\h(\Phi)\cong \A_\h(\Phi^u).$

One can prove  $R$-analogues of Corollary \ref{open1} and of Lemma
\ref{hphi}.
This way, one deduces that
the natural projection $\A_\h(\Phi^u)/\h\cdot\A_\h(\Phi^u)\iso R[x,y,z]$ is a graded
algebra isomorphism.
Thus, the algebra $\A_\h(\Phi^u)$
provides  a $\ct$-equivariant {\em flat} formal   deformation (where
$\hbar$ is the deformation parameter and where
the  $\ct$-equivariant structure comes from the grading)
 of   $\raa:=R[x,y,z],$
the latter being viewed as a Poisson $R$-algebra with an
$R$-bilinear Poisson bracket
arising from the polynomial $\phi^u:=(\Phi^u_1)^{\oper{ab}},$ cf.
\eqref{defbra}.

Recall next that to any  formal deformation-quantization of a commutative
algebra $A$ one can associate
a  Poisson bivector $\pi_\h\in\X^2A[[\h]]$ that represents the  Kontsevich's
class
of the deformation. The
Kontsevich correspondence is known to respect equivariant
structures arising from a   reductive group action by Poisson automorphisms.
Furthermore,  according to a result of Dolgushev \cite{Do}, the bivector
$\pi_\h$ gives a  {\em unimodular}  Poisson structure if and only if
the corresponding  deformation-quantization gives a flat family
of CY algebras. These results by Kontsevich and Dolgushev  can be easily
 generalized
to the setting of (flat) $R$-algebras.

Now, put $\raa_\h:=\raa[[\h]]$
and let $\pi_\h\in\X^2_R\raa_\h$ be a Poisson bivector
  that represents  Kontsevich's class of the deformation-quantization
of $\raa$ provided by the noncommutative $R[[\h]]$-algebra $\A_\h(\Phi^u)$.
 We know, by $R$-analogues of Corollary \ref{open1} and of Lemma
\ref{hphi}, that this deformation is indeed a  flat family
of CY $R$-algebras.
Therefore, we conclude that
the  $ R[[\h]]$-bilinear Poisson bracket $\{-,-\}$ on
$\raa_\h$ associated with the bivector $\pi_\h$ is
 unimodular. Moreover, since the
Kontsevich correspondence respects the $\ct$-equivariant
structure arising from the grading on $ \A_\h(\Phi^u)$,
resp. on $\raa_\h$,
we deduce that the  Poisson bracket associated
with the bivector $\pi_\h$ respects the grading
on the algebra $\raa_\h$, i.e. is such  that we have $\deg \{f,g\}
=\deg f+\deg g,$
for any homogeneous elements $f,g\in\raa_\h$ (where $\deg\h=0$ as before).

Next, one uses an $R$-analogue of Corollary \ref{alg_closed}(i)
to show
that there exists a formal series of the form
$\psi=\h\cdot\psi_1+\h^2\cdot\psi_2+\ldots,\
\psi_j\in\raa,$ such that, in $\X^2_R\raa[[\h]]$, one has
$\pi_\h=i_{\dr\psi}\ups$.
Here, $\ups\in\X^3_R\raa$ is
the  standard 3-vector  given
by formula \eqref{Vol}. Thus, $\deg\ups=-(a+b+c)=-d$.
It follows that each element $\psi_j\in\raa^{(d)}$ must be  homogeneous
 of degree $d$. It is also immediate from \eqref{defbra} that,
for the first term in the expansion of $\psi$, one has
\beq{phipsi}
\psi_1=(\Phi^u_1)^{\oper{ab}}.
\eeq

We introduce
$\raa\hb $, a commutative $\rkk$-algebra.
One may obviously view $\psi$ as an element of
$\raa\hb $.  Associated with
 this element,  there is a Poisson
$\rkk$-algebra $\raa_\psi$, cf. Definition \ref{aabb}. Clearly, we have
$\raa_\psi\cong \rkk\bigotimes_{R[[\h]]}\raa_\h,$ and the
Poisson bracket
on ${\mathcal{RA}}_{\psi}$ is nothing but the $\rkk$-bilinear
extension of the Poisson bracket  on the $\h$-algebra
$\raa_\h.$ Similarly, associated with the potential
$\Phi^u$, we have an $\rkk$-algebra
 $\A(\Phi^u):=\rkk\bigotimes_{R[[\h]]}\A_\h(\Phi^u)$.

At this point, one applies
Kontsevich's formality theorem \cite{Ko}, cf. also \cite{CVB}. The theorem yields
a  graded $\rkk$-algebra isomorphism,
\beq{konts}
\HH^\hd(\A(\Phi^u)\big)\cong \PH^\hd({\mathcal{RA}}_{\psi}).
\eeq

In degree zero, in particular,
we get  algebra isomorphisms
$\dis Z(\A(\Phi^u)\big)\cong \zz({\mathcal{RA}}_{\psi})
=\rkk[\psi].$
We deduce that the center of $\A(\Phi^u)$ is
generated by a degree $d$ homogeneous  element.
 Multiplying by a power of $\h$, we may assume without loss of generality
that this central element has the form
$1\o\Psi^u\in \rkk\bigotimes_{R[[\h]]}\A^{(d)}_\h(\Phi^u),$
where $\Psi^u\in\A^{(d)}_\h(\Phi^u)$ is such that $\Psi^u\ (\text{mod}\ \h)=\psi$.
Notice further that the $\h$-algebra $\A_\h(\Phi^u)$ has no $\h$-torsion
since $\Phi^u$  is a CY-potential, see  Proposition \ref{open}(i).
It follows that the  map
$\A_\h(\Phi^u)\to\A(\Phi^u),\, a\mto 1\o a,$
is injective and therefore $\Psi^u$ must be a central
element of the algebra  $\A_\h(\Phi^u)$.

By construction, the original potential $\Phi$ is
obtained
by specializing the homogeneous  potential $\Phi^u$ at $u=1$.
Thus we see that specializing the  central element $\Psi^u$ at $u=1$
one obtains a central
element $\Psi\in\A(\Phi)$, as required in the statement of  Proposition
\ref{cen}(i).
\end{proof}

\subsection{Proof of Proposition \ref{cen}(ii)-(iii)}\label{proof_cen}
Part (ii) is also an
easy consequence of the Kontsevich isomorphism \eqref{konts}. However, assuming
the statement of part (i) holds, one can give an alternate proof
of  part (ii) which does not involve formality theorem.
To this end, we exploit an adaptation of an argument used by
  Van den Bergh in the proof of \cite{VB2}, Theorem
4.1.

Recall that $\pi=i_{\dr\phi}\ups,$ cf. \eqref{Vol}.
First, we need the following corollary of Pichereau's results.

\begin{lem}\label{10}
The algebra $\PH^\hdot(\aa_\phi)$ is   a graded commutative
algebra with generators
$$\phi\in \PH^0(\aa_\phi),\quad \eu\in\PH^1(\aa_\phi),\quad \theta_1=i_{\dr f_1\!}\ups,\ldots,
\theta_{\mu-1\!}=i_{\dr f_{\mu-1}\!}\ups,\pi
\in \PH^2(\aa_\phi),\quad \ups\in\PH^3(\aa_\phi),
$$
and the following defining relations
\beq{relations}
\eu\cup \pi=\psi\cd \ups,\quad
\eu\cup \ups=\pi\cup \ups=0,\quad
\theta_i\cup \theta_j=
\theta_i\cup \pi=0,\en \forall i,j.
\eeq
\end{lem}

\begin{proof} For any polynomial $f\in \k[x,y,z]$, we have
$\eu\wedge i_{\dr f}\ups=\eu(f)\cdot\ups$.
Hence, we deduce $\eu\wedge i_{\dr\phi}\ups=d\cdot\phi\cdot\ups$.
Similarly, we get $\eu\wedge i_{\dr f_k}\ups=(\deg f_k)\cdot f_k\cdot\ups,$
for any $k=1,\ldots,\mu-1$. The statement readily follows from this
using the description of Poisson cohomology given in
Proposition  \ref{HPA}.
\end{proof}

Next,
we let
$\A_\h(\Phi)\supset\h\cdot\A_\h(\Phi)\supset\h^2\cdot\A_\h(\Phi)\supset\ldots,$
be the standard $\h$-adic filtration. The latter may be
extended in a unique way to a descending
$\Z$-filtration on the algebra $\A_\kk(\Phi)$ such that multiplication by
$\h\inv$ shifts the  filtration by one and such that
for the associated graded algebra, we have
$\gr\A_\kk(\Phi)=R[x,y,z][\h,\h\inv].$

The resulting associated graded Poisson bracket on $\gr\A_\kk(\Phi)$
is easily seen to be
the $\C[\h,\h\inv]$-bilinear extension of
the  Poisson bracket $\{-,-\}_\phi,$ on $\aa_\phi$,
where $\phi=(\Phi_1)^{\text{ab}}$.
In other words, we have a Poisson $\k[\h,\h\inv]$-algebra isomorphism
$\gr\A_\kk(\Phi)\cong\aa_\phi[\h,\h\inv].$

Associated with the above defined descending  filtration on
the   algebra $\A_\kk(\Phi)$, there is a standard
spectral sequence with  the first term,
cf. \cite{VB2}, page 224,
\beq{E1}
E_1=\PH^\hdot(\gr\A_\kk(\Phi))=\C[\h,\h\inv]\o\PH^\hdot(\aa_\phi)
\quad\Longrightarrow\quad\gr\HH^\hd(\A_\kk(\Phi)).
\eeq

Following an idea of Van den Bergh, we prove

\begin{lem}\label{degenerate}
Each of the elements from the
set of generators of the algebra $\PH^\hd(\aa_\phi)$
given in Lemma \ref{10} can be lifted to an element
in $\HH^\hd(\A_\kk(\Phi))$ in such a way that analogues of relations
\eqref{relations} hold for the lifted elements as well.
\end{lem}
\begin{proof}[Proof of Lemma] Set $\A_\kk=\A_\kk(\Phi).$
By Proposition \ref{cen}(i)
we have $\HH^0(\A_\kk )=\kk[\Psi]$. Furthermore,  the central element
$\Psi\in \A_\kk $ provides a lift of the element
$\phi\in \aa_\phi,$ due to equation \eqref{phipsi}.

To lift cohomology classes of degree 3, we compare
two duality isomorphisms provided by Proposition \ref{weinstein}(i)
and by \eqref{vdb}, respectively:
\begin{eqnarray*}
g:&\aa_\phi/\{\aa_\phi,\aa_\phi\}=\PH_0(\aa_\phi)\iso
\PH^3(\aa_\phi);\\
G:&\A_\kk /[\A_\kk,\A_\kk ]=\HH_0(\A_\kk )
\iso\HH^3(\A_\kk).
\end{eqnarray*}

Observe that any element  $f\in\aa_\phi/\{\aa_\phi,\aa_\phi\}$
can be trivially lifted to an element
$F\in\A_\kk /[\A_\kk ,\A_\kk ]$.
It follows easily that  any class of the form $g(f)\in \PH^3(\aa_\phi)$
admits a lift of the form $G(F)\in\HH^3(\A_\kk) $.
Further, let $\BB(F)\in \HH_1(\A_\kk )$
be the image of $F$ under the Connes differential
$\BB: \HH_0(\A_\kk )\to\HH_1(\A_\kk )$.
Then, one shows that
$G(\BB(F))\in \HH^2(\A_\kk ),$
the image of  $\BB(F)$ under the duality \eqref{vdb},
provides a lift of the class $i_{\dr f}\ups\in \PH^2(\aa_\phi)$.
In particular, each of the Poisson cohomology classes
$\pi=i_{\dr\phi}\ups,$ resp. $\theta_k=i_{\dr f_k}\ups, k=1\ldots,\mu-1,$
in $\PH^2(\aa_\phi)$,
see Lemma \ref{10},
has a lift $\Pi=G(\BB(\Psi)),$ resp. $\Theta_k=G(\BB(F_k)),$
in $\HH^2(\A_\kk ).$

Finally,  one lifts the class $\eu\in\PH^1(\aa_\phi)$ to
the corresponding Euler derivation $\Eu$ of the
{\em graded} algebra $\A_\kk $.

It follows from our construction that all of the relations from
\eqref{relations}, except possibly the first one,
 automatically hold for the lifted elements, by degree reasons.
Also the remaining relation holds for it is a
 a formal
consequence of \cite{Gi}, Theorem 3.4.3(i) and
the equation $\Eu(\Psi)=d\cdot\Psi.$
\end{proof}

According to the lemma, the  assignment sending our generators
of the algebra $\PH^\hd(\aa_\phi)$ to the corresponding
 generators
of the algebra $\HH^\hd(\A_\kk(\Phi))$
 can be extended to a well
defined graded
$\kk$-algebra map $\rho:\
\kk\o\PH^\hd(\aa_\phi)\to\HH^\hd(\A_\kk(\Phi))$.

We claim that the  map $\rho$ is an isomorphism.
To prove this, we exploit  \cite{VB2}, Lemma 5.2.
That lemma, combined with
our Lemma \ref{degenerate}, implies that the spectral sequence
in \eqref{E1} degenerates at  $E_1$. We deduce that, for the  filtration on
$\HH^\hd(\A_\kk(\Phi))$ coming from the spectral sequence,
one has
\beq{e1}
\gr\HH^\hd(\A_\kk(\Phi))\cong E_1=\C[\h,\h\inv]\o\PH^\hdot(\aa_\phi).
\eeq

Observe further that the lifts constructed
in Lemma  \ref{degenerate} are compatible with the filtrations involved.
Moreover, each term of the filtration is complete in the
$\h$-adic topology.
This, together with isomorphism \eqref{e1} immediately implies,
as explained at the top of page 224 in \cite{VB2},
that the map $\rho$ must be a bijection. That completes the proof of
part (ii) of Proposition \ref{cen}  and, at the same time, yields
part (iii), cf. Proposition \ref{HPA}(i).
\qed

\subsection{Proof of  Theorem \ref{one} and Theorem  \ref{three}}\label{1and3}

Part (i) of Theorem \ref{one}
follows directly from Corollary \ref{open1}(i) and Proposition
\ref{open}(ii).

Next, we prove  the existence of a central element in $\A(\Phi)$
from Theorem \ref{one}(ii) for
 {\em generic } potentials $\Phi\in F^{\leq d}\nat,$ where  $d=a+b+c$.
To this end,
one may
replace the ground field $\k$ by a larger field
and follow the strategy  of Van den Bergh, \cite{VB2}, \S5.
Thus,
we let our  ground field be of the form
$K(\!(\h)\!)$, for a certain field $K$.

We  assume (as we may) that the  coefficitients
in the expansion of $\Phi$ as a linear combination
of cyclic monomials in $x,y,z$  are algebraically independent over $\Q$.
Then, following \cite{VB2}, \S5, we may assume
that the potential has the form $\Phi=xyz-yxz+\sum_{j>0}\h^j\cdot\Phi_j,$
where $\Phi_j\in F^{\leq d}\nat$.
In such a case, Proposition \ref{cen}(i) insures the existence
of a  central element $\Psi\in \A^{\leq d}(\Phi),$ and we are done.

The proof of part (ii) of Theorem \ref{one} in the general case is based on
 a continuity argument. We
will  use the same notation concerning Rees algebras
as in the proof of Proposition \ref{cen}(i).

Thus, given a potential $\Phi=\Phi^{(d)}+\Phi^{(d-1)}+\ldots+\Phi^{(0)}$
of degree $\leq d$, we replace it by a degree $d$ homogeneous
potential
 $\Phi^u=\Phi^{(d)}+u\cdot\Phi^{(d-1)}+\ldots+u^d\cdot\Phi^{(0)}
\in \r F\nat$, where $\deg u=1$.
Further,
given  $\wt\Psi^u\in \r F^{(d)}$  let $\Psi^u\in \A(\Phi^u)$ denote
the image of $\wt\Psi$ under the projection
$\r F^{(d)}\onto\A^{(d)}(\Phi^u)$.
The condition that $\Psi^u\in \A(\Phi^u)$
be a central element of the algebra $\A(\Phi^u)$
amounts to the following 3 constraints on $ \wt\Psi^u,$
\beq{constr} v\cd\wt\Psi^u-\wt\Psi^u\cd v\in \pa_\Phi(\D^{(d+\deg v)}(\Phi^u)_1),\qquad
\forall v\in \{x,y,z\}.
\eeq
The commutator on the left is taken in the algebra $\r F,$
and $\pa_\Phi: \D^{(d)}(\Phi^u)_1\to\D^{(d)}(\Phi^u)_0=F$ stands for the
differential in the dg algebra $\D(\Phi^u)$.

Let $\Xi\sset \r  F^{(d)}\nat \times \mathbb P(\r F^{(d)})$
be the set of pairs $(\Phi^u,\k\cd \wt\Psi^u),$ where
$\Phi^u\in\r  F^{(d)}\nat$ is  a homogeneous  CY-potential
and the element $\wt\Psi^u$, generating
the line $\k\cd \wt\Psi^u\sset \r F^{(d)},$ satisfies \eqref{constr}.
According to  \eqref{rank}, for each
$r\geq 0$, the dimension of the vector
space $\pa_\Phi(\D^{(r)}(\Phi^u)_1)$ is a (finite) integer
independent of the choice of a CY-potential $\Phi^u\in \r  F^{(d)}\nat.$
It follows that the first projection
$\Xi\to \r  F^{(d)}\nat,\,(\Phi^u,\k\cd \wt\Psi^u)\mto \Phi^u,$ is
a proper morphism. The image of this morphism is dense
in $\r F^{(d)}\nat$ since we have already established the existence
of central elements in $\A(\Phi^{\leq d})$ for  generic potentials.
We conclude that the map $\Xi\to\r  F^{(d)}\nat$ is surjective,
and our claim follows by the specialization $u\mto1,\
\Phi^u\mto\Phi,$ and $\Psi^u\mto\Psi.$
\qed

\begin{proof}[Proof  of  Theorem  \ref{three}] Part (1) has been proved earlier,
in \S\ref{proof_one}.  To prove (2), we  repeat
the argument used in the proof of  Theorem  \ref{one}
in the case of generic potentials. This way, we see that the required
statement follows from the statement of Proposition \ref{cen}(iii) about
the center of the algebra $\A_\kk(\Phi)$.
\end{proof}

\subsection{Proof of Theorem \ref{four}}\label{pf4}
The statement of part (i) is a graded version of
the corresponding statement of  Theorem \ref{three}(1).
Thus, it follows from the latter theorem.

To prove part (ii), we may again
 reduce the statement to the case where
the ground field is  $\kk=\C(\!(\h)\!)$.
Furthermore, we may assume the
potential $\Phi$ to be  of the form \eqref{expand} and such that
 $(\Phi_1)^{\oper{ab}}\in\k[x,y,z]$ is a generic homogeneous  polynomial
of degree $d$.
Our assumptions on the triple $(a,b,c)$ insure that such a polynomial
 has an isolated singularity.
Thus, we are in a position
to apply  Proposition \ref{cen}(ii). The statement of
Theorem \ref{four}(ii)
then follows from that   proposition  and from
the corresponding results about Poisson
cohomology proved by Pichereau and
summarized in Proposition \ref{HPA}.

We now prove Theorem \ref{four}(iii). We keep the above setting, in
particular,
 we have $\kk$ as the base field.
Thus,  $A=\A_\kk(\Phi)$ is a  Calabi-Yau algebra and
we know that $\HH^1(A)=\kk[\Psi]\Eu,$
by  Proposition \ref{cen}(iii).

Let $\vol\in \HH_3(A)$ denote the image of $1\in Z(A)=\HH^0(A)$ under
 the duality isomorphism \eqref{vdb}.
Then, the duality gives a
$\kk[\Psi]$-module isomorphism
$\HH^1(A)\iso\HH_2(A)$ that sends $\Eu\in \HH^1(A)$ to $i_{\Eu}\vol\in\HH_2(A)$.
Therefore, using the equation
$\Eu(\Psi)=d\cdot\Psi$ and standard calculus identities
in the framework of Hochschild cohomology, cf. \cite{Lo}, \S 4.1, we
 compute (where dot denotes cup-product on Hochschild cohomology),
\begin{eqnarray*}
&\BB(\Psi^k\cd i_{\Eu}\vol)&=
\BB\ccirc i_{\Eu}(\Psi^k\cd \vol)=
(\BB\ccirc i_{\Eu}+ i_{\Eu}\ccirc \BB)(\Psi^k\cd \vol)\\
&&=L_\Eu(\Psi^k\cd \vol)=k\cd\Psi^{k-1}\cd\Eu(\Psi)\cd \vol
+\Psi^k\cd L_\Eu\vol\\
&&=kd\cd\Psi^k\cd \vol+d\cd\Psi^k\cd\vol=
(k+1)d\cd\Psi^k\cd \vol.
\end{eqnarray*}

Since $(k+1)d\neq 0$ for any $k=0,1,\ldots,$
from the calculation above we deduce that the Connes differential
gives a bijection $\BB: \HH_2(A)\iso \HH_3(A).$
By duality, this implies that the BV-differential
yields a bijection $\De: \HH^1(A)\iso\HH^0(A)$.
That proves one of the two isomorphisms
of Theorem \ref{four}(iii).

To prove the other isomorphism, we
observe  that $A$ is a
nonnegatively  graded algebra
with degree zero component
equal to $\kk.$
Hence,
by \cite{EG}, Lemma 3.6.1, we get an exact
sequence of Hochschild homology
\beq{deltaexact1}
0\to\kk\to \HH_0(A)\stackrel{\BB}\too \HH_1(A)\stackrel{\BB}\too \HH_2(A)\stackrel{\BB}\too
\HH_3(A)\to 0.
\eeq
Applying duality \eqref{vdb}, we obtain an exact sequence of  Hochschild cohomology
\beq{deltaexact}0\to\kk\cd\ups\to \HH^3({A})\stackrel{\Delta}\too \HH^2({A})
\stackrel{\Delta}\too \HH^1({A})\stackrel{\Delta}\too
\HH^0({A})\to 0.
\eeq

We have shown earlier that the last map $\De$ on the right in this exact sequence
is a bijection. This forces the first map  $\De$ on the left to be
a surjection, and we are done.\qed

\begin{rem}\label{sur} There are also
Poisson cohomology counterparts of exact sequences
\eqref{deltaexact1}-\eqref{deltaexact}.  The  counterpart of
\eqref{deltaexact1}
follows, using Cartan's
homotopy formula $L_\eu=\dr\ccirc i_\eu+i_\eu\ccirc\dr,$ from
the fact that the operator $L_\eu$ acts on $\Om^j\aa_\phi$ with positive
weights,
for any $j\neq 0.$ The  analogue of \eqref{deltaexact} can be deduced
from this by
duality,
cf. Proposition \ref{weinstein}.

Further, an explicit description of the group $\PH^2(\aa_\phi)$
given  by Pichereau \cite{P} shows that the map
$\delta: \PH^3(\aa_\phi)\to \PH^2(\aa_\phi)$,
equivalently, the map $\dr: \PH_0(\aa_\phi)\to\PH_1(\aa_\phi)$, is surjective as well.
This, combined with  spectral sequence
\eqref{e1}, may be used to obtain an alternate proof
of  Theorem \ref{four}(iii).
\end{rem}

\subsection{Sketch of proof of  Theorem \ref{two}}\label{sketch2}

We  begin with part (i).
First of all we introduce a space of  deformation
parameters similar to the one used
 in the proof of Theorem \ref{delpezzoform}.
Specifically, let $S_\A$ be  the space of tuples
$(t,c,P,Q,R)$. We have
$\dim S_\A=(p-1)+(q-1)+(r-1)+2=\mu,$ by \eqref{dimJ}.

For each $s=(t,c,P,Q,R)\in S_\A$ we let $A_s:=\A(\Phi^{t,c}_\pqr)$
be the corresponding algebra. This is a filtered algebra,
with an associated graded algebra $\gr A_s.$
Hence there is an induced ascending filtration
$\HH^\hd_{\leq m}(A_s)$ on  Hochschild cohomology, resp. homology, groups of
$A_s$.
Proving  Theorem \ref{two}(i) amounts to showing that
there exists a subset  $U\sset S_\A$, of sufficiently general parameters,
  such that for any $s\in U,$ the  Kodaira-Spencer  map
induces an isomorphism
\beq{ks1}
\KS_s:\ T_sS_\A\iso \HH^2_{\leq 0}(A_s),
\qquad\forall s\in U \  (\sset S_\A).
\eeq

To this end, we first
 use the classification result from Theorem \ref{three}(i).
The theorem implies that for any choice of
subset $F^\circ\nat\sset F\nat$,
of generic potentials  in the
sense of Definition \ref{generic_def},
the set $U:=\{(t,c,P,Q,R)\in S_\A\;|\;\Phi^{t,c}_\pqr\in F^\circ\nat\}$
is nonempty and, moreover, it is a subset of generic parameters in
 $S_\A$,
 in the
sense of Definition \ref{generic_def} again.

We have the following
diagram, cf.
\eqref{kappadeform},
\beq{KS1}
\xymatrix{
T_sS_\A\ar[d]_<>(0.5){\KS_s}\ar[rr]^<>(0.5){pr}&&(A_s)\nat
\ar@{->>}[d]^<>(0.5){\mathsf{B}}\\
\HH^2(A_s)\ &&\ \HH_1(A_s)
\ar@{=}[ll]^<>(0.5){\eqref{vdb}}\\
}
\eeq

In this diagram,
the map $pr$ is the tautological projection that sends
a variation of the potential, viewed as an element of
$\k\langle x,y,z\rangle\nat$, to its image
in $(A_s)\nat$. Observe further that the  isomorphism
\eqref{vdb} at the bottom of
the diagram gives a bijection between $\HH_\idot^{\leq d}(A_s)$ and
$\HH^{d-\hd}_{\leq 0}(A_s)$. Furthermore, Proposition
\ref{KS_prop} insures that diagram \eqref{KS1} commutes.

In order to prove
\eqref{ks1} for an algebra $A_s$ associated
with a  potential $\Phi=\Phi^{t,c}_\pqr$ with
{\em generic}
coefficients,   we may (and will) assume  that
our base field is $\kk=\k(\!(\h)\!)$ and that our potential
has the form \eqref{expand}. We put $\phi:=(\Phi_1)^{\text{ab}},$
cf. \eqref{defbra} and
let  $\aa_\phi$ be the corresponding  Poisson algebra.

There is an analogue of  diagram \eqref{KS1} for
the Poisson algebra $\aa_\phi$ instead of the algebra $A_s.$ Furthermore, there is
a spectral sequence like \eqref{e1} for each of
the Hochschild (co)homology groups in \eqref{KS1}.
Its $E_1$-term is the corresponding
Poisson  (co)homology group in the Poisson analogue of \eqref{KS1}.

First of all, applying Proposition \ref{cen}(ii) we get
$\dim \HH^2_{\leq 0}(A_s)=\dim \PH^2_{\leq 0}(\aa_\phi)$.
Now, for any homogeneous element $f$ and $k\geq 0$, we have
$\deg(\phi^k\cdot i_{\dr f}\ups)=kd+\deg f-(a+b+c)=\deg f+ (k-1)d.$ Therefore,
using Proposition
\ref{weinstein} and the notation of Proposition \ref{HPA}(ii), we find that the elements
$\dr f_1,\ldots,\dr f_{\mu-1},\,\dr\phi,$ form a $\C$-basis of the vector
space $\PH_1^{\leq d}(\aa_\phi).$ Thus, we deduce
\beq{bound}
\dim \HH^2_{\leq 0}(A_s)=\dim \PH^2_{\leq 0}(\aa_\phi)=
\dim \PH_1^{\leq d}(\aa_\phi)=\mu=\dim S_\A.
\eeq

Thus, to complete the proof of part (i) it suffices to show that the
map \eqref{ks1} is surjective.
 From diagram \eqref{KS1},
we see that this would follow provided we prove the surjectivity of the
composite map
${\BB}\ccirc pr: T_sS_\A\to \HH_1^{\leq d}(A_s).$
Using the spectral sequence in \eqref{e1} we reduce the
latter statement to proving surjectivity of a similar
map $T_sS_\A \to\PH_1^{\leq d}(\aa_\phi),$
for Poisson algebras.  But this is clear since there are obvious
elements in $f_j\in S_\A=\C^2\times S_p\times S_q\times S_r,$
cf. \S\ref{ellip}, proof of Theorem \ref{delpezzoform}, such that the 1-forms
$\dr f_1,\ldots,\dr f_{\mu-1},\dr \phi,$  give a basis of the vector
space $\PH_1^{\leq d}(\aa_\phi).$

The proof of Theorem \ref{three}(ii) proceeds in a similar way.
We omit the details.
\qed

\medskip
{\large
\section{Appendix: computer calculation}\label{magma_sec}}
\subsection{}
In the $E_6$ case the relations in the algebra
$\A(\Phi^{t,c}_\pqr)$ take the following form
\begin{align*}
&x  y-q   y  x- t   z^2+c_1   z+c_2,\\
&y z-q  z  y- t  x^2+a_1  x+a_2,\\
&z  x-q  x  z- t  y^2+b_1  y+b_2
\end{align*}

The corresponding central element $\Psi$ was computed by Eric Rains
using MAGMA. It reads
\begin{align*}
& t(q+1)(t(t^3+1)y^3 + (q^3 -   t^3)  y
z   x - q(t^3+1)
  z   y   x + t(q^3 -   t^3)  z^3)\\
&\en\en-
    t(q^2 +
    q  t^3 + q + 2  t^3 + 1) b_1  y^2\\
&\en\en +
 (q  t^3-q^2)
  a_1  y  z +
     t^3(q+1)  b_1  z x +
(q^3+ q  t^3)  a_1  z y \\
&\en\en+q(q+1)  t^3 c_1
  y   x+t(2 q  t^3+
    t^3-q^4 - q^3 - q^2) c_1  z^2 \\
&\en\en-  ((q^3  t  + 2  q^2  t+q  t )  a_2
  + q^2  a_1^2 +   q  t^2  b_1  c_1 )  x \\
&\en\en- t((q^3  b_2  + 2  q^2   + q  t^3 +2q +t^3 +1 )b_2 +  q  t
 a_1  c_1     -   t^2  b_1^2)  y \\
&\en\en-t((q^4+ 2  q^3+ 2  q^2- q t^3 + q  -t^3)  c_2
           +   q  t^2
 c_1^2 +  q  t  a_1  b_1    )  z
\end{align*}

We refer to \cite{R} for more complicated formulas in the $E_7$
and $E_8$ cases.

\begin{rem}
We were informed by the referee that such formulas were also 
obtained by a computer calculation in D. Stephenson's thesis. 
\end{rem}


{\small
\begin{thebibliography}{9999}
\bibitem[A]{A} M.
Artin, {\em
Some problems on three-dimensional graded domains.}
 Representation theory and algebraic geometry. 1--19,
London Math. Soc. Lecture Note Ser., 238, Cambridge Univ. Press, 1997.


\bibitem[AS]{AS} \bysame, W. Schelter,
{\em Graded algebras of global dimension $3$.}  Adv. in Math.  66  (1987),   171--216.
\bibitem[ATV]{ATV} \bysame, J. Tate, M. Van den Bergh,
{\em  Some algebras associated to automorphisms of elliptic curves.}
  The Grothendieck Festschrift, Vol. I,  33--85, Progr. Math., 86,
  Birkh\"auser Boston,   1990.


\bibitem[AV]{AV} \bysame, M. Van den Bergh,
{\em Twisted homogeneous coordinate rings.}  J. Algebra  133  (1990),
249--271.
\bibitem[BT]{BT} R.  Berger, R.  Taillefer,
{\em Poincare-Birkhoff-Witt Deformations of Calabi-Yau Algebras.}
J. Noncomm. Geom.
{\bf 1} (2007), 241-270.
{\tt{arXiv:math.RT/0610112}}.


\bibitem[Bo]{Bo} R. Bocklandt,
{\em Graded Calabi Yau Algebras of dimension 3.}
  J. Pure Appl. Algebra  212  (2008),  14--32.
{\tt{arXiv:math.RA/0603558}}.

\bibitem[B]{B} E.
Brieskorn, {\em Singular elements of semi-simple algebraic groups.}
  Actes du Congr\`es International des Math\'e-maticiens (Nice, 1970), Tome
 {\bf  2},  pp. 279--284. Gauthier-Villars, Paris, 1971.

\bibitem[BSV]{BSV} L.Le Bruyn, P. Smith, M. Van den Bergh, {\em
Central extensions of three-dimensional Artin-Schelter regular algebras.}
Math. Z. {\bf  222} (1996), 171--212.

\bibitem[Br]{Br} J.-L. Brylinski,
{\em A differential complex for Poisson manifolds.}
  J. Differential Geom.  {\bf 28}  (1988),  93--114.

\bibitem[CVB]{CVB} D.  Calaque, M. Van den Bergh,
{\em Hochschild cohomology and Atiyah classes}.
{\tt{arXiv:0708.2725}}.

\bibitem[C1]{C1} T. Cassidy, {\em Global dimension 4 extensions of
  Artin-Schelter regular algebras.}
  J. Algebra {\bf 220} (1999), 225--254.

\bibitem[C2]{C2} \bysame, {\em
Central extensions of Stephenson's algebras.}
Comm. Algebra {\bf 31} (2003), 1615--1632.

\bibitem[CK]{CK} D.
Chan, R. Kulkarni, {\em
del Pezzo orders on projective surfaces.}
Adv. Math. {\bf  173} (2003), 144--177.



\bibitem[CBEG]{CBEG} W. Crawley-Boevey, P. Etingof, V.
 Ginzburg, {\em
Noncommutative Geometry and Quiver algebras.} Adv. in Math.
{\bf 209} (2007), 274-336.
{\tt{arXiv:math.AG/0502301}}.

\bibitem[Do]{Do} V. Dolgushev,
{\em The Van den Bergh duality and the modular symmetry of a Poisson
variety.}
 Selecta Math. (N.S.)  14  (2009),  199--228. 
 {\tt{arXiv:math/0612288}}.
 \bibitem[D]{D} {\em
S\'eminaire sur les Singularit\'es des Surfaces.}
1976--1977.
Lect. Notes in Math., {\bf 777}. Springer-Verlag, 1980.

\bibitem[Ei]{Ei} D. Eisenbud, {\em
Homological algebra on a complete intersection, with an application to
group representations.}
  Trans. Amer. Math. Soc. {\bf  260} (1980), 35--64.
 \bibitem[EG]{EG} P. Etingof, V. Ginzburg,
{\em  Noncommutative complete intersections and matrix
integrals.} Pure Appl. Math. Quart. {\bf 3} (2007), 107-151.
{\tt{arXiv:math/0603272}}.

\bibitem[EOR]{EOR} \bysame, A. Oblomkov, E. Rains,
{\em Generalized double affine Hecke algebras of rank 1 and quantized Del Pezzo surfaces}
 Adv. Math.  212  (2007), 749--796. {\tt{arXiv:math/0406480}}.

\bibitem[Gi]{Gi} V. Ginzburg, {\em Calabi-Yau algebras},
{\tt{arXiv:math.AG/0612139}}.

\bibitem[GK]{GK}  \bysame, D. Kaledin, {\em
Poisson deformations of symplectic quotient singularities.}
Adv. Math. {\bf 186} (2004), 1--57.


\bibitem[KST]{KST} H.
 Kajiura, K. Saito, A. Takahashi, {\em
Matrix Factorizations and Representations of Quivers II: type ADE case.}
 Adv. Math.  211  (2007), 327--362. {\tt{arXiv:math/0511155}}.
\bibitem[Ka]{Ka} D.
Kaledin, {\em On the coordinate ring of a projective Poisson scheme.}
 Math. Res. Lett.  13  (2006),  99--107.
\bibitem[K1]{Ko} M. Kontsevich,
{\em Deformation quantization of Poisson manifolds.}
  Lett. Math. Phys.  66  (2003),  no. 3, 157--216.
{\tt{arXiv:q-alg/9709040}}.
\bibitem[K2]{Ko2} M. Kontsevich,
{\em  Deformation quantization of algebraic varieties}.
  Lett. Math. Phys.  56 (2001), 271-294.

\bibitem[LR]{LR} A. Lago, A. Rodicio,
{\em Generalized Koszul complexes and Hochschild (co-)homology of
 complete intersections.}
  Invent. Math. {\bf  107}  (1992),  433--446.
\bibitem[LPP]{LPP} R.
Laza, G. Pfister, D. Popescu, {\em
Maximal Cohen-Macaulay modules over the cone of an elliptic curve.}
J. Alg. {\bf 253} (2002), 209--236.
\bibitem[Le]{Le} T.
Levasseur, {\em Some properties of noncommutative regular graded rings.}
 Glasgow Math. J. {\bf  34}  (1992),  277--300.




\bibitem[Lo]{Lo} J.L. Loday, {\em   Cyclic homology.}
 Grundlehren der Mathematischen Wissenschaften, 301.
  Springer-Verlag,  1998.
\bibitem[Ma]{Ma} N. Marconnet, {\em Homologies of cubic Artin-Schelter
 regular algebras.}
 J. Alg. {\bf  278}  (2004),   {638-~665.}


\bibitem[NVB]{NVB} K. de Naeghel, M. Van den Bergh,
{\em Ideal classes of three-dimensional Sklyanin algebras.}
  J. Algebra  {\bf 276}  (2004),   515--551.

\bibitem[Or1]{Or1} D. Orlov,
{\em Derived categories of coherent sheaves and triangulated categories
of singularities.}  Proc. Steklov Inst. Math.  2004, (246), 227--248.
{\tt{ arXiv:math/0302304}}.

\bibitem[Or2]{Or2}\bysame, {\em
 Triangulated categories of singularities and D-branes in Landau-Ginzburg models.}
  Tr. Mat. Inst. Steklova  246  (2004),  Algebr. Geom. Metody, Svyazi i
 Prilozh., 240--262;
translation in  Proc. Steklov Inst. Math.  2004,  (246), 227--248.
{\tt{arXiv:math/0503632}}.

\bibitem[P]{P} A. Pichereau, {\em  Poisson (co)homology and isolated singularities.}
 J. Algebra  {\bf 299}  (2006),   747--777. cf. also
 C. R. Math. Acad. Sci. Paris {\bf 340}  (2005),  151--154.
{\tt{arXiv:math/0511201}}.
\bibitem[R]{R} E. Rains, Computer calculation available at\
{\text{http://www-math.mit.edu/{\~\,}\!etingof/delpezzocenter}.}

\bibitem[Sa]{Sa} K. Saito,
{\em  Regular system of weights and associated singularities.}
Complex analytic singularities,  479--526, Adv. Stud. Pure Math., 8, North-Holland, Amsterdam, 1987.

\bibitem[St1]{St1} D. R. Stephenson, {\em  Algebras Associated to Elliptic Curves.} Trans.
Amer. Math. Soc. {\bf 349} (1997), 2317-2340.

\bibitem[St2]{St2} D. R. Stephenson, {\em Artin-Schelter Regular Algebras of Global Dimension Three,} PhD Thesis, 
University Microfilms International, 1994.

\bibitem[St3]{St3} D. R. Stephenson, {\em Artin-Schelter Regular Algebras of Global Dimension Three,} 
Journal of Algebra 183, (1996), 55-73. 

 \bibitem[VB1]{VB1} M. Van den Bergh,
{\em  A relation between Hochschild homology and cohomology for
Gorenstein rings.}
Proc. Amer. Math. Soc. {\bf 126} (1998),  1345--1348;  130  (2002),
 2809--2810.
\bibitem[VB2]{VB2} \bysame, {\em Noncommutative homology of some three-dimensional
 quantum spaces.}  $K$-Theory {\bf  8}  (1994),   213--230.

\bibitem[VB3]{VB3} \bysame, {\em
 Blowing up of non-commutative smooth surfaces.}
  Mem. Amer. Math. Soc. {\bf  154}  (2001),  no. 734.

\bibitem[Xu]{Xu} P.
Xu, {\em
Gerstenhaber algebras and BV-algebras in Poisson geometry.}
Comm. Math. Phys. {\bf 200} (1999), 545--560.
\end{thebibliography}
}
\end{document}